\documentclass[10pt,oneside,leqno]{amsart}
\usepackage{latexsym,amsmath,amsthm,amssymb,enumerate,mathrsfs,graphicx}
\usepackage{color}
\usepackage[matrix,arrow]{xy}
\newtheorem{theorem}{Theorem}
\newtheorem{proposition}{Proposition}

\theoremstyle{definition}
\newtheorem{remark}{Remark}

\begin{document}

\title[Six dimensional Lie algebras]{Cohomological properties of unimodular six dimensional solvable Lie algebras}
\author{Maura Macr\`i}
\address{Dipartimento di Matematica G. Peano \\ Universit\`a di Torino\\
Via Carlo Alberto 10\\
10123 Torino\\ Italy} \email{maura.macri@unito.it}


\begin{abstract} In the present paper we study six dimensional solvable Lie algebras with special emphasis on those admitting a symplectic structure. We list all the symplectic structures that they admit and we compute  their Betti numbers finding some properties about the codimension of the nilradical.
Next,  we consider the conjecture of Guan about step of nilpotency of a symplectic solvmanifold finding that it is true for all six dimensional unimodular solvable Lie algebras. 
Finally, we consider some cohomologies for symplectic manifolds introduced by Tseng and Yau in the context of symplectic Hogde theory and we use them to determine some six dimensional solvmanifolds for which the Hard Lefschetz property holds.
 \end{abstract}

\maketitle

\section*{Introduction}
A solvmanifold $M=G/\Gamma$ is a compact homogeneous space of a solvable Lie group $G$, i.e. a compact quotient of a solvable Lie group $G$ by a lattice $\Gamma$. 
A special class of solvmanifolds, called nilmanifolds, was introduced by Malcev \cite{Ma} and corresponds to the particular case when $G$ is a nilpotent Lie group.

\smallskip

Both classes of manifolds have been particularly important for producing examples of compact symplectic manifolds which do not admit any K\"ahler structure. 
Hence the study of the topology  of solvmanifolds 
(and in particular their cohomological properties) is particularly interesting, especially when they are endowed with a symplectic structure. In this context, the Hard Lefschetz property, formality and symplectic Hodge theory play an important role (see for instance \cite{FOT}). 

\smallskip

Hence one needs to compute the de Rham cohomology of a solvmanifold, which, in some situations can be done using invariant differential forms, i.e. by the  Chevalley-Eilenberg cohomology $H^* ({\mathfrak g})$ of ${\mathfrak g}$.  This is the case if  the \textit {Mostow condition} holds, namely, the algebraic closures ${\mathcal A} (\mbox{Ad}_G (G))$ and  ${\mathcal A} (\mbox{Ad}_G (\Gamma))$ are equal,  \cite{mostow2} and \cite[Corollary 7.29]{rag}.
Special instances are provided by nilmanifolds \cite{nomizu} and completely solvable Lie groups $G$ \cite{hattori}, i.e.,  the adjoint representation ${\mbox {ad}}_X: {\mathfrak g} \to {\mathfrak g}$ have only real eigenvalues for all $X \in {\mathfrak g}$.

\smallskip

Unlike nilpotent Lie groups there is no simple criterion to understand whether a solvable Lie group $G$ admits a lattice, but a necessary condition is that $G$ is \textit{unimodular},  i.e., for all $X\in \mathfrak{g}, \; \mbox{tr} \,  \mbox {ad}_{X}=0$, where $ \mathfrak{g} $ is the Lie algebra of $G$.

\smallskip

Nilpotent Lie algebras and solvable Lie algebras have been classified up to dimension 5 (see for instance \cite{bock}). In dimension 6, the number of possible solvable Lie algebras is very large. A  complete classification can be obtained by using the results in  the papers \cite{bock,muba,shaba}.  In this paper we improve some classifications, considering only  the case of six dimensional solvable unimodular Lie algebras.

\smallskip

Solvmanifolds up to dimension six admitting an invariant symplectic structure (invariant means that it comes from a symplectic form on the Lie algebra) were studied by Bock \cite{bock}. In particular, he considered the conditions of being cohomologically symplectic, formality and Hard Lefschetz property.

 \smallskip
 
 In the present paper we study six dimensional  unimodular solvable (non-nilpotent) Lie algebras with special emphasis on those admitting a symplectic structure. Six dimensional nilpotent Lie algebras admitting symplectic structures were classified in \cite{salamon}.

In Section 1  we list all the symplectic structures  that six dimensional  unimodular solvable (non-nilpotent) Lie algebras admit  (Table 3)  and we consider the conjecture of Guan \cite{guan} about steps of a symplectic solvmanifold \cite{guan}, namely  that if  a solvmanifold  $G/\Gamma$ admits a symplectic structure  then $G$ is  at most 3-step solvable.  We find that this is true for a six dimensional Lie group whose Lie algebra is unimodular, indeed this holds for all six dimensional unimodular solvable Lie algebra, also those which do not admit symplectic structures, (see Proposition~\ref{steps-guan}).

 In Section 2  we compute  their Betti numbers (i.e., the dimensions of their Chevalley-Eilenberg cohomology) finding some properties about the codimension of the nilradical. Recall that the \textit{nilradical} of a Lie algebra $\mathfrak{g}$ is its largest  nilpotent ideal.

Namely, in Section~\ref{betti_numbers} we prove as main result

\begin{theorem}\label{betti}
Let $ \mathfrak{g} $ be a six dimensional unimodular, solvable, non-nilpotent Lie algebra 
\begin{itemize}
\item if it admits a  symplectic structure, then it has positive, non zero, second Betti number, i.e., $ b_{2}(\mathfrak{g})>0 $.
\item if  $ b_{1}(\mathfrak{g})=1 $, then its nilradical has codimension 1 and  $ b_{2}(\mathfrak{g})=0 $ if and only if $ b_{3}(\mathfrak{g})=0 $.
\item if its  nilradical has codimension greater then 1, then $ b_{1}(\mathfrak{g})\geq 2 $ and $ b_{2}(\mathfrak{g})=1 $ if and only if $ b_{3}(\mathfrak{g})=0 $.
\end{itemize}
\end{theorem}

The Betti numbers of the 6-dimensional Lie algebras with 5-dimensional nilradical  were also computed by M. Freibert and F. F. Schulte-Hengesbach \cite{FS}.

\smallskip

\smallskip

Finally, we consider the Hard Lefschetz property and some cohomologies for symplectic manifolds introduced by Tseng and Yau \cite{tsengyau} in the context of symplectic Hogde theory. In particular, we show that these cohomologies can be computed using invariant forms, provided this is the case for the Rham cohomology (see Theorem~\ref{cinesi} in Section~\ref{hard_lefschetz}).
We apply this result, together with the list of symplectic structures on solvable Lie algebras (Table 3), to show that some solvmanifolds satisfy the Hard Lefschetz property (Theorem~\ref{HL}).
\smallskip 

In the Appendix we include two Tables. In Table 2 we list all the solvable, non-nilpotent unimodular six dimensional Lie algebras, indicating the differential of the generators of the dual algebra. This list is based on the classifications given in \cite{bock}, \cite{muba} and \cite{shaba} and has the same notation for the Lie algebras. 
 In Table 3 we list all the symplectic structures on six dimentional solvable unimodular Lie algebras.

\medskip

\textcolor{red}{ \noindent {\em{Acknowledgements}}. I would like to thank Antonio Otal and Anna Fino for the precious help given in the last version of the paper. }

\medskip
\newpage

\section{Symplectic Structures}\label{sympl_str}
Let $\mathfrak{g}$ be a real Lie algebra of dimension $2n$.   We recall that a  \textit{symplectic structure} on $\mathfrak{g}$ is a closed $2$-form $\omega$ in $ \bigwedge \mathfrak{g}^* $ such that $ \omega^n\neq 0 $.
Let $\mathfrak{g}$ be a six dimensional real solvable unimodular Lie algebra and let $ \{X_1,...,X_6\} $ be an ordered basis of $\mathfrak{g}$, then a $2$-form $\omega$ is associated in a natural way to a matrix $ M=(\omega_{ij}) \in \mathcal{M}(6,\mathbb{R}) $, where $ \omega_{ij}:= \omega (X_i,X_j) $, and $ \omega^n\neq 0 \Leftrightarrow \det M \neq 0 $. 

\medskip

By direct computation we prove

\begin{theorem}
The six dimensional real solvable, non-nilpotent unimodular Lie algebras admitting a  symplectic structure are the following:
\\{\small
$\mathfrak{g}_{6.3}^{0,-1}=(-26,-36,0,-46,56,0), \newline
\mathfrak{g}_{6.10}^{0,0}=(-26,-36,0,-56,46,0), \newline
\mathfrak{g}_{6.13}^{\frac{1}{2},-1,0}=(-23+\frac{1}{2}.16,-\frac{1}{2}.26,36,-46,0,0), \newline
\mathfrak{g}_{6.13}^{-1,\frac{1}{2},0}=(-23+\frac{1}{2}.16,26,-\frac{1}{2}.36,-46,0,0), \newline
\mathfrak{g}_{6.15}^{-1}=(-23,-26,36,-26-46,-36+56,0), \newline
\mathfrak{g}_{6.18}^{-1,-1}=(-23,26,-36,-36-46,56,0), \newline
\mathfrak{g}_{6.21}^{0}=(-23,0,-26,-46,56,0), \newline
\mathfrak{g}_{6.23}^{0,0, \varepsilon}=(-23-\varepsilon.56,0,-26,-36,0,0),  \quad \varepsilon \neq 0\newline
\mathfrak{g}_{6.29}^{0,0,\varepsilon}=(-23-\varepsilon.56,0,0,-36,-46,0), \newline
\mathfrak{g}_{6.36}^{0,0}=(-23,0,-26,56,-46,0), \newline
\mathfrak{g}_{6.38}^{0}=(-23,36,-26,-26+56,-36-46,0), \newline
\mathfrak{g}_{6.54}^{0,-1}=(-35-16,-45+26,-36,46,0,0), \newline
\mathfrak{g}_{6.70}^{0,0}=(-35+26,-45-16,46,-36,0,0), \newline
\mathfrak{g}_{6.78}=(-25+16,-45,-24-36-46,-46,56,0), \newline
\mathfrak{g}_{6.118}^{0,\pm 1,-1}=(-25+16,15+26,\mp 45-36,\pm 35-46,0,0), \newline
\mathfrak{n}_{6.84}^{\pm 1}=(-45,-15-36,-14+26\mp 56,56,-46,0), \newline
\mathfrak{g}_{5.7}^{p,-p,-1}\oplus \mathbb{R}=(-15,-p.25,p.35,45,0,0), \newline
\mathfrak{g}_{5.8}^{-1}\oplus \mathbb{R}=(-25,0,-35,45,0,0), \newline
\mathfrak{g}_{5.14}^{0}\oplus \mathbb{R}=(-25,0,-45,35,0,0), \newline
\mathfrak{g}_{5.17}^{0,0,r}\oplus \mathbb{R}=(-25,15,-r.45,r.35,0,0), \newline
\mathfrak{g}_{5.17}^{p,-p,\pm 1}\oplus \mathbb{R}=(-p.15-25,15-p.25,p.35\mp 45,\pm 35+p.45,0,0), \newline
\mathfrak{g}_{5.17}^{0,0,\pm 1}\oplus \mathbb{R}=(-25,15,\mp 45,\pm 35,0,0), \newline
\mathfrak{g}_{5.18}^{0}\oplus \mathbb{R}=(-25-35,15-45,-45,35,0,0), \newline
\mathfrak{g}_{5.19}^{-2,2}\oplus \mathbb{R}=(-23+15,-25,+2.35,-2.45,0,0), \newline
\mathfrak{g}_{5.19}^{-\frac{1}{2},-1}\oplus \mathbb{R}=(-23-\frac{1}{2}.15,-25,-\frac{1}{2}.35,45,0,0), \newline
\mathfrak{g}_{3.4}^{-1}\oplus 3\mathbb{R}=(-13,23,0,0,0,0), \newline
\mathfrak{g}_{3.5}^{0}\oplus 3\mathbb{R}=(-23,13,0,0,0,0), \newline
\mathfrak{g}_{3.1}\oplus \mathfrak{g}_{3.4}^{-1}=(-23,0,0,-46,56,0), \newline
\mathfrak{g}_{3.1}\oplus \mathfrak{g}_{3.5}^{0}=(-23,0,0,-56,46,0), \newline
\mathfrak{g}_{3.4}^{-1}\oplus \mathfrak{g}_{3.4}^{-1}=(-13,23,0,-46,56,0), \newline
\mathfrak{g}_{3.4}^{-1}\oplus \mathfrak{g}_{3.5}^{0}=(-13,23,0,-56,46,0), \newline
\mathfrak{g}_{3.5}^{0}\oplus \mathfrak{g}_{3.5}^{0}=(-23,13,0,-56,46,0)$},\newline
where the parameters $\varepsilon,  p$ and $r$ are real numbers.
\end{theorem}

To explain this notation, for example \,{\small $\mathfrak{g}_{5.17}^{0,0,r}\oplus \mathbb{R}=(-25,15,-r.45,r.35,0,0)$} means that there is a basis $(\alpha_1, \dots , \alpha_6)$ of  the dual of the Lie algebra $\mathfrak{g}_{5.17}^{0,0,r}\oplus \mathbb{R}$ such that 
$ d\alpha_1= - \alpha_{25}, d\alpha_2=\alpha_{15}, d\alpha_3=-r\alpha_{45}, d\alpha_4=r\alpha_{35}, d\alpha_5=0, d\alpha_6=0$, where by $\alpha_{ij}$ we denote  $\alpha_i \wedge \alpha_j$.

\begin{proof}
To construct the symplectic form we take the generic element $\omega \in \ker d \subset \bigwedge^2 (\mathfrak{g}^*)$ and we impose it to be not degenerate, that is $ \omega^3\neq 0 $.
\\ With this direct computation we can see that the six dimentional solvable unimodular Lie algebras not listed below have always $ \omega^3 =0 $ for every $\omega \in \ker d \subset \bigwedge^2 (\mathfrak{g}^*)$. In Table 3 (Appendix) we list the symplectic structures admitted.
\end{proof}

Describing nilmanifolds and solvmanifolds with symplectic structure became important after the work of Thurston, \cite{thurston}. For this reason in \cite{guan}, Guan 
studied properties about the steps of nilmanifolds, showing that if a nilmanifold $G/\Gamma$ admits a symplectic structure then $G$ has to be at most two step as a solvable Lie group. He also 
conjectured that the Lie group of a solvmanifold admitting a symplectic structure is at most 3-step solvable. 

We show by direct computation that this is true for all six dimensional unimodular solvable Lie algebra, also for  those  which do not admit any symplectic structure.

\begin{proposition}\label{steps-guan}
Every six dimensional unimodular, solvable, non-nilpotent Lie algebra $ \mathfrak{g} $ is 2 or 3-step solvable, in particular
\begin{itemize}
\item if its nilradical has codimension 1, it is 3-step solvable unless it is almost abelian, or $\mathfrak{g}$ is one of the following Lie algebras: \newline $ \mathfrak{g}_{6.14}^{a,0}, \quad  \mathfrak{g}_{6.17}, \quad   \mathfrak{g}_{6.18}^{0,0}, \quad  \mathfrak{g}_{6.20}, \quad  \mathfrak{g}_{6.21}^{0}, \quad  \mathfrak{g}_{6.23}^{0,0,\varepsilon}, \quad  \mathfrak{g}_{6.25}^{-1,0}, \quad  \mathfrak{g}_{6.29}^{0,0,\varepsilon}, \newline  \mathfrak{g}_{6.36}^{0,0}, \quad  \mathfrak{g}_{6.54}^{0,-1}, \quad  \mathfrak{g}_{6.63}, \quad  \mathfrak{g}_{6.65}^{0,0}, \quad  \mathfrak{g}_{6.70}^{0,0}, \quad  \mathfrak{g}_{6.88}^{0,0,0}$.
\item if its nilradical has codimension greater then 1, it is 2-step solvable unless $\mathfrak{g}$ is one of the following Lie algebras: \newline $\mathfrak{g}_{6.129}, \quad  \mathfrak{g}_{6.135}, \quad  \mathfrak{g}_{5.19}\oplus \mathbb{R}, \quad \mathfrak{g}_{5.20}\oplus \mathbb{R}, \quad \mathfrak{g}_{5.23}\oplus \mathbb{R}, \quad \mathfrak{g}_{5.25}\oplus \mathbb{R}, \quad \mathfrak{g}_{5.26}\oplus \mathbb{R}, \quad \mathfrak{g}_{5.28}\oplus \mathbb{R}, \quad \mathfrak{g}_{5.30}\oplus \mathbb{R}, \quad \mathfrak{g}_{4.8}\oplus 2\mathbb{R}, \quad \mathfrak{g}_{4.9}\oplus 2\mathbb{R}$.
\end{itemize}
\end{proposition}

\bigskip 
\bigskip 
\bigskip

\section{Betti numbers of 6-dimensional unimodular solvable\\ non-nilpotent Lie algebras} \label{betti_numbers}
In this Section we compute the second and third Betti number of six dimensional solvable Lie algebras. The interest in determining solvable Lie algebras with the property that $b_2(\mathfrak{g})= b_3(\mathfrak{g})$ comes from a class of manifolds endowed with a closed 3 form, called \textit{String geometry}, considered in \cite{mad-sw2}.
Strong geometry is an important example of connection between mathematics and physics, in particular multi-moment maps are used in string theory and one-dimensional quantum mechanics, \cite{mad-sw2}.

\medskip 

Let $ M $ be a manifold, then $ (M,\gamma) $ is a \textit{Strong geometry} if $ \gamma $ is a closed 3-form on $ M $. 
Suppose there is a Lie group $ G $ that acts on $ M $ preserving $ \gamma $, then we denote by $\mathit{P}_{\mathfrak{g}}$ the kernel of the map $ \bigwedge^2 \mathfrak{g} \rightarrow \mathfrak{g}$ induced by the Lie bracket of $ \mathfrak{g} $. \\A \textit{Multi-moment map} is an equivariant map $ \nu : M  \rightarrow \mathit{P}_{\mathfrak{g}}^*$ such that $d \langle \nu, p \rangle = i_p \gamma$, for any $p \in\mathit{P}_{\mathfrak{g}}$, (where $i_p$ denotes the interior product) \cite{mad-sw1}.
\\ We refer to \cite{mad-sw1} and \cite{mad-sw2} for details on strong geometry. In particular Madsen and Swann \cite{mad-sw2} proved that 
if $ b_2(\mathfrak{g})= b_3(\mathfrak{g})=0 $, then there exists a multi-moment map for the action of $ G $ on the manifold $M$.
Because of this result  they listed the Lie algebras with trivial second and third Betti numbers, up to dimension five.
We add to their classification the Betti numbers of 6-dimensional solvable, non-nilpotent unimodular Lie algebras.

\begin{remark} Every Lie algebra $\mathfrak{g}$ whose Lie group is solvable has $ b_1(\mathfrak{g})>0 $, \cite{bock}.
\end{remark}

Next we list  6-dimensional unimodular, solvable, non-nilpotent Lie algebra $ \mathfrak{g}$ together with their first, second and third Betti number. The Betti numbers of the 6-dimensional Lie algebras with 5-dimensional nilradical  were also computed by M. Freibert and F. F. Schulte-Hengesbach \cite{FS}.

\medskip 
Looking at this list and comparing with Table 3 yields Theorem \ref{betti}.

\bigskip

\centerline{{\sc{Table 1}}: {\small{Betti numbers of 6 dimensional unimodular, solvable, non-nilpotent Lie algebras \footnote{In Table 1 we impose conditions which become at every step more restrictive. It is therefore implicit that the previous conditions hold only when the more restrictive ones are not satisfied.}}}}

\begin{center}
{\small
\begin{tabular}{|c|l @{}  p{1.4cm}|l @{} p{4cm}|l @{} p{4cm}|}
\hline
$ \mathfrak{g} $ & $ b_{1} $ & & $ b_{2} $ & & $ b_{3} $ & \\
\hline
$ \mathfrak{g}_{6.1} $ & 1 & &
0 \;& if $ a\neq -1, \: b\neq -1,\: b\neq -a,$ & 0 \; & if $ a\neq -1,\: b\neq -1,\: b\neq -a, $ \\
&&&& $ c\neq -a, \: c+b\neq -1,\: c\neq -b,  $ && $ c\neq -a, \: c+b\neq -1,\: c\neq -b, $ \\
&&&& $a+b\neq -1,\:a+c\neq -1$ && $a+b\neq -1,\:a+c\neq -1 $ \\
 & & & 1\; & if $ a = -1 $,  or if $ b=-a $, & 2 \; & if $ a = -1 $,  or if $ b=-a $, \\
 & & & & or if $ b=-c$, or if $a+b \neq -1$ & &  or if $  b=-c$, or if $a+b \neq -1$\\
  & & & 2\; & if $ b = -1 $, & 4 \; & if $ b = -1 $, \\
 & & & & or if $c =-a $ or if $ c =-1-a $, &  & or if $c =- a $ or if $ c=-1-a $,\\
 & & & & or if $a=-1$ and $b=1$, & &  or if $a=-1$ and $b=1$,\\
 & & & & or if $a=-1$ and $b+c=-1$,  & &  or if $a=-1$ and $b+c=-1$,\\  
 & & & & or if $ b=c= -a$, & & or if $ b=c= -a$,\\
 & & & & or if $b=-c=\pm a$, & &  or if $b=-c=\pm a$,\\
 & & & & or if $ b=-c=\pm (1+a)$, & & or if $ b=-c=\pm (1+a)$,\\
 & & & & or if $b=c=-1-a$ & &  or if $b=c=-1-a$ \\
 & & & 3\; & if $a=-\frac{1}{2}$ and $b =-c = \pm \frac{1}{2} $ & 6 \; & if $a=-\frac{1}{2}$ and $b =-c = \pm \frac{1}{2} $ \\
& & & & or if $a=b=c-\frac{1}{2}$ & &  or if  $a=b=c-\frac{1}{2}$\\
& & & 4\; & if $a=-b =c = \frac{1}{2} $ & 5 \; & if  $a=-b =c = \frac{1}{2} $ \\
 \hline
 $ \mathfrak{g}_{6.2} $ & 1 & if $ a \neq 0 $ &0\; &if $ a\neq 0,\: c\neq -1, \:e \neq -c$, &0\; &if $ a\neq 0,\: c\neq -1, \:e \neq -c, $ \\
 &&&& $c-e \neq \pm 1,\: c+e \neq 1 $, && $c-e \neq \pm 1,\: c+e \neq 1 $,\\
 & 2 & if $ a =0 $ &1 &if $ c=-1 $, &2 &if $ a=0 $ or if $c=-1$  \\
 &&&&or if $ e=-c $, or if  $ e = c+1 $,  &&or if $ e=-c $, or if  $ e = c+1 $, \\
 &&&& or if $ e=\pm (1-c) $,  && or if $ e=\pm (1-c) $, \\
 & & & 2 & if $ a =0 $ or if $ c=-e=\pm \frac{1}{2} $ & 4 & if $ c=-e=\pm \frac{1}{2} $ \\
  &&&& or if $e=-1$,  && or if $e=-1$,  \\
&&&& or if $c=-e=-1$ && or if $c=-e=-1$ \\
 \hline
 $ \mathfrak{g}_{6.3} $ &1&if $ a\neq -1 $&0& if $ a\neq -1,\frac{1}{2} $&0& if $ a\neq -1,\frac{1}{2} $\\
 &2&if $ a= -1 $ &1&if $ a=\frac{1}{2} $&2&if $ a=\frac{1}{2} $\\
 &&&3&if $ a=-1$&4&if $ a=-1 $\\
 \hline
 $ \mathfrak{g}_{6.4} $ &1&&0&&0 & \\
 \hline
 $ \mathfrak{g}_{6.6} $ &1&if $ a\neq -\frac{1}{2} $&0&if $ a\neq -1,-\frac{1}{2} $&0 &if $ a\neq -1,-\frac{1}{2} $ \\
 &2&if $ a= -\frac{1}{2} $ &1& if $ a=-1 $ & 2& if $ a=-1 $ \\
 &&&2& if $ a=-\frac{1}{2} $ & 2& if $ a=-\frac{1}{2} $ \\
 \hline

\end{tabular}
}
\end{center}

\begin{center}
{\small
\begin{tabular}{|c|l @{}  p{1.5cm}|l @{} p{4cm}|l @{} p{4cm}|}
\hline
$ \mathfrak{g} $ & $ b_{1} $ & & $ b_{2} $ & & $ b_{3} $ & \\
\hline
$ \mathfrak{g}_{6.7} $ &1&&0&&0 & \\
\hline
$ \mathfrak{g}_{6.8} $ &1&&0&if $ a+b\neq 0,\;a+c\neq 0,$&0 &if $ a+b\neq 0,\;a+c\neq 0,$ \\
&&&& $b+c\neq 0,\;p\neq 0 $ && $b+c\neq 0,\;p\neq 0$\\
&&&1& if $ a+b=0 $ ,&2& if $ a+b=0 $,\\
&&&& or if $ b+c=0 $,&& or if $ b+c=0 $,\\
&&&& or if $ p=0 $ && or if $ p=0 $\\
&&&2& if  $ a=-b=c $,&4& if $ a=-b=c $ ,\\
&&&& or if $ a+c=0 $ &&or if $ a+c=0 $\\
\hline
$ \mathfrak{g}_{6.9} $ & 1 & if $ b\neq 0 $&0& if $ bp\neq 0, \; a+b\neq 0 $ &0& if $ bp\neq 0, \; a+b\neq 0 $\\
&2& if $ b=0 $ &1& if $ p=0 $ or if $ a+b=0 $ &2& if $ bp=0 $ or if $ a+b=0 $\\
&&&2& if $ b=0 $ &&\\
\hline
$ \mathfrak{g}_{6.10} $ & 1 & if $ a\neq 0 $&0& if $ a\neq 0$ &0& if $ a\neq 0$\\
&2& if $ a=0 $ &3& if $ a=0 $&4& if $ a=0 $\\
\hline
$ \mathfrak{g}_{6.11} $ & 1\; &&0\;& if $ pq\neq 0$ &0\;& if $ pq\neq 0$\\
&&&1& if $ pq=0 $&2& if $ pq=0 $\\
\hline
$ \mathfrak{g}_{6.12} $ &1&&0&&0 & \\
\hline
$ \mathfrak{g}_{6.13} $ &1& if $ bh\neq 0 $&0& if $ a\neq -1, \;b\neq -1,  $&0& if $ a\neq -1, \;b\neq -1, $\\
&&&& $a+b\neq 0, \;2a+b\neq 0,  $ && $ a+b\neq 0, \;2a+b\neq 0, $\\
&&&& $ a+2b\neq 0, \;a+2b+1\neq 0,$ && $ a+2b\neq 0, \;a+2b+1\neq 0,  $ \\
&&&& $ b+2a+1\neq 0   $ && $ b+2a+1\neq 0   $\\
&2& if $ b=0 $ &1& if $ a=-1 $ or if $b=-1$  &2 & if $ a=-1 $ or if $b=-1$  \\
&&or if $ h=0 $&& or if $a+b=0$ && or if $a+b=0$\\
&&&&  or if $a+2b=0,-1$ &&  or if $a+2b=0,-1$ \\
&&&& or if $b+2a=0,-1$ &&  or if $b+2a=0,-1$ \\
&&&2& if $ a=-1 $ and $ b=2 $&4 & if $ a=-1 $ and $ b=2$ \\
&&&& or if $ b =-1 $ and $ a=  2$&& or if $ b =-1 $ and $ a= 2$\\
&&&& or if $ a =\frac{1}{3} $ and $ b=-\frac{2}{3} $ && or if $ a =\frac{1}{3} $ and $ b=-\frac{2}{3} $\\
&&&& or if $ a= -\frac{2}{3}$ and $ b= \frac{1}{3}$  && or if $ a= -\frac{2}{3}$ and $ b= \frac{1}{3}$\\
&&&& or if $ a=b=-1,-\frac{1}{3} $ && or if $ a=b=-1,-\frac{1}{3} $\\
&&&3& if $ a=\frac{1}{2}$ and $ b=-1 $&4& if  $ a=\frac{1}{2}$ and $ b=-1 $ \\
&&&& or if $a=-1$ and $b=0,\frac{1}{2}$ &&  or if $a=-1$ and $b=0,\frac{1}{2}$ \\
&&&& or if $a=-b=\pm 1$ &6& if $a=-b=\pm 1$ \\
\hline
$ \mathfrak{g}_{6.14} $&1& if $ a\neq -\frac{1}{3} $&0& if $ a\neq -1, -\frac{2}{3}, -\frac{1}{3}, \frac{1}{3}, \frac{2}{3} $&0& if $ a\neq -1, -\frac{2}{3}, \frac{1}{3}, \frac{2}{3} $\\
&2& if $ a= -\frac{1}{3} $&1& if $ a= -1, -\frac{2}{3}, -\frac{1}{3}, \frac{1}{3}, \frac{2}{3} $&2& if $ a\neq -1, -\frac{2}{3}, \frac{1}{3}, \frac{2}{3} $\\
\hline
$ \mathfrak{g}_{6.15} $ &1&&2&&4 & \\
\hline 
$ \mathfrak{g}_{6.17} $ &2&&2&&1 & \\
\hline
$ \mathfrak{g}_{6.18} $ &1&if $ a\neq 0 $&0& if $ a\neq 0,-\frac{1}{2},-1,-2,-3 $&0 & if $ a\neq -\frac{1}{2},-1,-2,-3 $ \\
&2& if $ a=0 $&1 & if $ a= 0,-\frac{1}{2},-2,-3 $&2& if $ a= -\frac{1}{2},-2,-3 $\\
&&&2& if $ a=-1 $&4&if $ a=-1 $\\
\hline 
$ \mathfrak{g}_{6.19} $ &1&&0&&0 & \\
\hline 
$ \mathfrak{g}_{6.20} $ &2&&1&&0 & \\
\hline 
$ \mathfrak{g}_{6.21} $ &1& if $ a\neq 0 $&0& if $ a\neq 0,-\frac{1}{3},-1 $&0 &if $ a\neq 0,-\frac{1}{3},-1 $ \\
&2& if $ a=0 $ &1& if $ a=-\frac{1}{3},-1 $&2& if $ a=-\frac{1}{3},-1 $\\
&&&3& if $ a=0$&4& if $ a=0$\\
\hline 
$ \mathfrak{g}_{6.22} $ &1&&0&&0 & \\
\hline
$ \mathfrak{g}_{6.23} $ &1& if $ a\neq 0 $&0&if $ a\neq 0 $&0 & if $ a\neq 0 $\\
&3& if $ a=0 $&5&if $ a=0 $&6 &if $ a=0 $ \\
\hline
$ \mathfrak{g}_{6.25} $ &1& if $ b\neq 0,-1 $&0&if $ b\neq 0,-1,-\frac{1}{2},1 $&0 & if $ b\neq 0,-1,-\frac{1}{2},1$\\
&2& if $ b=0,-1 $&1&if $ b=-\frac{1}{2},1 $&2 &if $ b= 0,-1,-\frac{1}{2},1  $ \\
&&&2& if $ b=0,-1 $ &&\\
\hline 

\end{tabular}
}
\end{center}

\begin{center}
{\small
\begin{tabular}{|c|l @{}  p{2cm}|l @{} p{4cm}|l @{} p{4cm}|}
\hline
$ \mathfrak{g} $ & $ b_{1} $ & & $ b_{2} $ & & $ b_{3} $ & \\
\hline
$ \mathfrak{g}_{6.26} $ &2&&2&&2 & \\
\hline
$ \mathfrak{g}_{6.27} $ &1&&1&&2 & \\
\hline
$ \mathfrak{g}_{6.28} $ &1&&0&&0 & \\
\hline
$ \mathfrak{g}_{6.29} $ &1& if $ b\neq 0 $&2& if $ b \neq 0$&2 & if $ b \neq 0$ \\
&3& if $ b=0 $&5&  if $ b=0 $ and $ \varepsilon \neq 0 $&6&  if $ b=0 $ and $ \varepsilon \neq  0 $\\
&&&6&  if $ b=0 $ and $ \varepsilon = 0 $&8& if $ b=0 $ and $ \varepsilon=0 $\\
\hline
$ \mathfrak{g}_{6.32} $ &1\;& if $ a\neq -\frac{h}{2},-\frac{h}{6} $&0\;& if $ a \neq  0,-\frac{h}{2},-\frac{h}{6}$&0\; & if $ a \neq 0$ \\
&2\;& if $ a= -\frac{h}{2},-\frac{h}{6} $&1&  if $ a= 0,-\frac{h}{2},-\frac{h}{6} $ &2&  if $ a=0 $ \\
\hline
$ \mathfrak{g}_{6.33} $ &1& if $ a\neq 0 $&0& if $ a\neq 0 $&0 & if $ a\neq 0 $ \\
&3& if $ a= 0 $&3& if $ a=0 $&1 & if $ a=0 $ \\
\hline
$ \mathfrak{g}_{6.34} $ &1& if $ a\neq 0 $&0& if $ a\neq 0 $&0 & if $ a\neq 0 $ \\
&3& if $ a= 0 $&3& if $ a=0 $&1 & if $ a=0 $ \\
\hline

$ \mathfrak{g}_{6.35} $ &1& if $ ab\neq 0 $&0& if $ c\neq 0, a\neq 0,-2b, $&0 & if $ c\neq 0, a\neq -2b, $  \\
&&&& $ b\neq 0,-2a $&& $ b\neq -2a $\\
&2& if $ a= 0 $ &1& if $ a=0 $ or if $ b= 0 $ or if $ c=0 $ &2 & if  $ c=0 $ or if $ a=-2b $  \\
&&or $ b= 0 $&& or if $ a=-2b $ or if $ b=-2a $&& or if $ b=-2a $\\
\hline
$ \mathfrak{g}_{6.36} $ &1& if $ a\neq 0 $&0& if $ a\neq 0 $&0 & if $ a\neq 0 $ \\
&2& if $ a= 0 $&3& if $ a=0 $&4 & if $ a=0 $ \\
\hline
$ \mathfrak{g}_{6.37} $ &1&&0& if $ a\neq 0 $&0 & if $ a\neq 0 $ \\
&&&1& if $ a=0 $&2 & if $ a=0 $ \\
\hline
$ \mathfrak{g}_{6.38} $ &1&&2&&4 & \\
\hline
$ \mathfrak{g}_{6.39} $ &1& if $ h\neq 0 $&0& if $ h\neq 0,-\frac{1}{2},-1,-2,-3 $&0 & if $ h\neq -\frac{1}{2},-1,-2,-3 $ \\
&2& if $ h=0 $&1& if $ h=0,-\frac{1}{2},-1,-2,-3  $&2 & if $ h=-\frac{1}{2},-1,-2,-3  $ \\
\hline
$ \mathfrak{g}_{6.40} $ &1&&0&&0& \\
\hline
$ \mathfrak{g}_{6.41} $ &1&&1&&2& \\
\hline
$ \mathfrak{g}_{6.42} $ &1&&0&&0& \\
\hline
$ \mathfrak{g}_{6.44} $ &1&&0&&0& \\
\hline
$ \mathfrak{g}_{6.47} $ &2&&1&&0& \\
\hline
$ \mathfrak{g}_{6.54} $ &1& if $ l\neq -\frac{1}{2},$&0& if $ l\neq 0,-\frac{1}{2},-1,-2,-\frac{3}{2},-\frac{2}{3} $&0 & if $ l\neq 0,-1,-\frac{3}{2},-\frac{2}{3}  $ \\
&&$ -1,-2  $&1& if $ l=0,-\frac{1}{2},-2,-\frac{3}{2},-\frac{2}{3}  $&2 & if $ l=0,-\frac{3}{2},-\frac{2}{3} $ \\
&2& if $ l=-\frac{1}{2},$&3& if $ l=-1 $&4& if $ l=-1 $\\
&& $ -1,-2  $&&&&\\
\hline
$ \mathfrak{g}_{6.55} $ &1&&0&&0& \\
\hline
$ \mathfrak{g}_{6.56} $ &1&&0&&0& \\
\hline
$ \mathfrak{g}_{6.57} $ &1&&0&&0& \\
\hline
$ \mathfrak{g}_{6.61} $ &1&&0&&0& \\
\hline
$ \mathfrak{g}_{6.63} $ &2&&2&&2& \\
\hline
$ \mathfrak{g}_{6.65} $ &1& if $ l\neq 0 $&0& if $ l\neq 0 $&0&  if $ l\neq 0 $\\
&3& if $ l=0 $&5& if $ l=0 $&6&  if $ l=0 $\\
\hline
$ \mathfrak{g}_{6.70} $ &1& if $ p\neq 0 $&0& if $ p\neq 0 $&0&  if $ p\neq 0 $\\
&2& if $ p=0 $&3& if $ p=0 $&4&  if $ p=0 $\\
\hline
$ \mathfrak{g}_{6.71} $ &1&&0&&0& \\
\hline
$ \mathfrak{g}_{6.76} $ &1&&1&&1& \\
\hline
$ \mathfrak{g}_{6.78} $ &1&&1&&2& \\
\hline
$ \mathfrak{g}_{6.83} $ &1& if $ l\neq 0 $&1& if $ l\neq 0 $&2&  if $ l\neq 0 $\\
&3& if $ l=0 $&5& if $ l=0 $&6&  if $ l=0 $\\
\hline
$ \mathfrak{g}_{6.84} $ &2&&2&&2& \\
\hline
$ \mathfrak{g}_{6.88} $ &1&if $ \mu\neq 0 $ or &1& if $ \mu\nu\neq 0 $&2& $ \mu\nu\neq 0 $ \\
&&$ \nu\neq 0 $ &3& if $ \mu\neq 0 $ and $ \nu=0 $&6& if $ \mu\neq 0 $ and $ \nu=0 $\\
&5& if $ \mu=0 $ and && or if $ \mu=0 $ and $ \nu\neq 0 $&& or if $ \mu=0 $ and $ \nu\neq 0 $\\
&&$ \nu=0 $&9& if $ \mu=0 $ and $ \nu=0 $&10\;& if $ \mu=0 $ and $ \nu=0 $\\
\hline

\end{tabular}
}
\end{center}

\begin{center}
{\small
\begin{tabular}{|c|l @{}  p{2cm}|l @{} p{4cm}|l @{} p{4cm}|}
\hline
$ \mathfrak{g} $ & $ b_{1} $ & & $ b_{2} $ & & $ b_{3} $ & \\
\hline
$ \mathfrak{g}_{6.89} $ &1&if  $ s\nu\neq 0 $&1& if $ s\nu\neq 0 $&2&if $ s\nu\neq 0 $ \\
&2& if $ s\neq 0 $&3& if $ s\neq 0 $ and $ \nu=0 $&& or if $ s\neq 0 $ and $ \nu=0 $\\
&& or $ \nu\neq 0 $ && or if $ s=0 $ and $ \nu\neq 0 $&& or if $ s=0 $ and $ \nu\neq 0 $\\
&5& if $ \nu=0 $ &9& if $ s=0 $ and $ \nu=0 $&10& if $ s=0 $ and $ \nu=0 $\\
&& and $ s=0 $ &&&&\\
\hline
$ \mathfrak{g}_{6.90} $ &1& if $ \nu\neq 0 $&1& if $ \nu\neq 0 $&2&\\
&3& if $ \nu=0 $&3& if $ \nu=0 $&&\\
\hline
$ \mathfrak{g}_{6.91} $ &1&&1&&2& \\
\hline
$ \mathfrak{g}_{6.92} $ &1\;&if  $ \mu\nu\neq 0 $&3\;& if $ \mu\nu\neq 0 $&4\;&if $ \mu\nu\neq 0 $ \\
&2& if $ \mu\neq 0 $&5& if $ \mu\neq 0 $ and $ \nu=0 $&6& if $ \mu\neq 0 $ and $ \nu=0 $\\
&& or if $ \nu\neq 0 $ && or if $\mu=0 $ and $ \nu\neq 0 $&& or if $ \mu=0 $ and $ \nu\neq 0 $\\
&5& if $ \nu=0 $ &9& if $ \mu=0 $ and $ \nu=0 $&10\;& if $ \mu=0 $ and $ \nu=0 $\\
&& and $ \mu=0 $ &&&&\\
\hline
$ \mathfrak{g}_{6.92}^{*} $ &1&&3&&6& \\
\hline
$ \mathfrak{g}_{6.93} $ &1& if $ \nu\neq 0 $&1& if $ \nu\neq 0 $&2& if $ \nu\neq 0 $\\
&3& if $ \nu=0 $&3& if $ \nu=0 $&2& if $ \nu=0 $\\
\hline
$ \mathfrak{g}_{6.94} $ &1&&1&&2& \\
\hline
$ \mathfrak{g}_{6.101} $ &2&&1& if $ a\neq -2 $ or $ b\neq -1 $ &0& if $ a\neq -2 $ or $ b\neq -1 $\\
&&&2& if $ a=-2 $ and $ b=-1 $&1& if $ a=-2 $ and $ b=-1 $\\
\hline
$ \mathfrak{g}_{6.102} $ &2&&1&&0& \\
\hline
$ \mathfrak{g}_{6.105} $ &2&&1&&0& \\
\hline
$ \mathfrak{g}_{6.107} $ &2&&1&&0& \\
\hline
$ \mathfrak{g}_{6.113} $ &2&&1& if $ a\neq 0 $ or $ b\neq -1 $ &0& if $ a\neq 0 $ or $ b\neq -1 $\\
&&&3& if $ a=0 $ and $ b=-1 $&2& if $ a=0 $ and $ b=-1 $\\
\hline
$ \mathfrak{g}_{6.114} $ &2&&2& if $ a\neq \pm 2 $ &2& if $ a\neq \pm 2 $ \\
&&&3& if $ a=\pm 2 $&3& if $ a=\pm 2 $\\
\hline
$ \mathfrak{g}_{6.115} $ &2&&1&&0& \\
\hline
$ \mathfrak{g}_{6.116} $ &2&&1&&0& \\
\hline
$ \mathfrak{g}_{6.118} $ &2&&1& if $ b\neq \pm 1 $ &0& if $ b\neq \pm 1 $ \\
&&&3& if $ b=\pm 1 $&4& if $ b=\pm 1 $\\
\hline
$ \mathfrak{g}_{6.120} $ &2&&2&&2& \\
\hline
$ \mathfrak{g}_{6.121} $ &2&&2&&2& \\
\hline
$ \mathfrak{g}_{6.129} $ &2&&1&&0& \\
\hline
$ \mathfrak{g}_{6.135} $ &2&&1&&0& \\
\hline
$ \mathfrak{n}_{6.83} $ &1&&1&&2& \\
\hline
$ \mathfrak{n}_{6.84} $ &1&&1&&1& \\
\hline
$ \mathfrak{n}_{6.96} $ &1& if $ b \neq 0 $&1& if $ b \neq 0 $&2& \\
&3& if $ b=0 $&3& if $ b=0 $&&\\
\hline
\hline 
$ \mathfrak{g}_{5.7}\oplus \mathbb{R} $&2\;&&1\;& if $ r\neq -1 $&0\;& if $ r\neq -1 $\\
&&&3&  if $ r=-1 $ and $ q\neq -1 $&4&  if $ r=-1 $ and $ q\neq -1 $\\
&&&5&  if $ r=-1 $ and $ q=-1 $&8&  if $ r=-1 $ and $ q=-1 $\\
\hline 
$ \mathfrak{g}_{5.8}\oplus \mathbb{R} $&3&&5&&6&\\
\hline 
$ \mathfrak{g}_{5.9}\oplus \mathbb{R} $&2\;& if $ p\neq 0 $&1\;& if $ p\neq 0,-1 $&0\;& if $ p\neq 0,-1 $\\
&3& if $ p=0 $&3&  if $ p=0,-1 $&2&  if $ p=0 $\\
&&&&&4&  if $ p=-1 $\\
\hline 
$ \mathfrak{g}_{5.11}\oplus \mathbb{R} $&2&&1&&0&\\
\hline
$ \mathfrak{g}_{5.13}\oplus \mathbb{R} $&2&&1& if $ q\neq 0 $&0& if $ q\neq 0 $\\
&&&3& if $ q=0 $&4& if $ q=0 $\\
\hline
$ \mathfrak{g}_{5.14}\oplus \mathbb{R} $&2&&5&&6&\\
\hline

\end{tabular}
}
\end{center}

\begin{center}
{\small
\begin{tabular}{|c|l @{}  p{2cm}|l @{} p{4cm}|l @{} p{4cm}|}
\hline
$ \mathfrak{g} $ & $ b_{1} $ & & $ b_{2} $ & & $ b_{3} $ & \\
\hline
$ \mathfrak{g}_{5.15}\oplus \mathbb{R} $&2&&3&&4&\\
\hline 
$ \mathfrak{g}_{5.16}\oplus \mathbb{R} $&2&&1&&0&\\
\hline
$ \mathfrak{g}_{5.17}\oplus \mathbb{R} $&2&&1& if $ p\neq 0 $ and $ r \neq \pm 1 $&0& if $ p\neq 0 $ and $ r \neq \pm 1 $ \\
&&&3& if $ p=0 $ and $ r\neq \pm 1 $&4& if $ p=0 $ and $ r\neq \pm 1 $ \\
&&&& or if $ p\neq 0 $ and $ r=\pm 1 $&& or if $ p\neq 0 $ and $ r=\pm 1 $\\
&&&5&  if $ p=0 $ and $ r=\pm 1 $&8& if $ p=0 $ and $ r=\pm 1 $ \\
\hline
$ \mathfrak{g}_{5.18}\oplus \mathbb{R} $&2&&3&&4&\\
\hline
$ \mathfrak{g}_{5.19}\oplus \mathbb{R} $&2& if $ p\neq 0 $&1& if $ p\neq 0,-\frac{1}{2},-2 $&0& if $ p\neq 0,-\frac{1}{2},-2 $\\
&3& if $ p=0 $&3& if $ p=0,-\frac{1}{2},-2 $&2& if $ p=0 $\\
&&&&&4& if $ p=-\frac{1}{2},-2 $\\
\hline 
$ \mathfrak{g}_{5.20}\oplus \mathbb{R} $&3&&3&&3&\\
\hline
$ \mathfrak{g}_{5.23}\oplus \mathbb{R} $&2&&1&&0&\\
\hline
$ \mathfrak{g}_{5.25}\oplus \mathbb{R} $&2&&1&&0&\\
\hline
$ \mathfrak{g}_{5.26}\oplus \mathbb{R} $&3&&3&&2&\\
\hline
$ \mathfrak{g}_{5.28}\oplus \mathbb{R} $&2&&1&&0&\\
\hline
$ \mathfrak{g}_{5.30}\oplus \mathbb{R} $&2&&1&&0&\\
\hline
$ \mathfrak{g}_{5.33}\oplus \mathbb{R} $&3&&3&&2&\\
\hline
$ \mathfrak{g}_{5.35}\oplus \mathbb{R} $&3&&3&&1&\\
\hline
\hline 
$ \mathfrak{g}_{4.2}\oplus 2\mathbb{R} $&3&&3&&2&\\
\hline
$ \mathfrak{g}_{4.5}\oplus 2\mathbb{R} $&3&&3&&2&\\
\hline
$ \mathfrak{g}_{4.6}\oplus 2\mathbb{R} $&3&&3&&2&\\
\hline
$ \mathfrak{g}_{4.8}\oplus 2\mathbb{R} $&3&&3&&2&\\
\hline
$ \mathfrak{g}_{4.9}\oplus 2\mathbb{R} $&3&&3&&2&\\
\hline
\hline 
$ \mathfrak{g}_{3.4}\oplus 3\mathbb{R} $&4&&7&&8&\\
\hline
$ \mathfrak{g}_{3.5}\oplus 3\mathbb{R} $&4&&7&&8&\\
\hline
$ \mathfrak{g}_{3.1}\oplus \mathfrak{g}_{3.4} $&3&&5&&6&\\
\hline
$ \mathfrak{g}_{3.1}\oplus \mathfrak{g}_{3.5} $&3&&5&&6&\\
\hline
$ \mathfrak{g}_{3.4}\oplus \mathfrak{g}_{3.4} $&2&&3&&4&\\
\hline
$ \mathfrak{g}_{3.4}\oplus \mathfrak{g}_{3.5} $&2&&3&&4&\\
\hline
$ \mathfrak{g}_{3.5}\oplus \mathfrak{g}_{3.5} $&2&&3&&4&\\
\hline

\end{tabular}
}
\end{center}
\bigskip 
\bigskip 
\bigskip

\bigskip 
\section{Hard Lefschetz property of 6 dimensional unimodular non-nilpotent solvmanifolds}\label{hard_lefschetz}
L.S. Tseng and S.T. Yau introduced some classes of finite dimensional cohomologies for symplectic manifolds \cite{tsengyau}. These cohomology classes depend on the symplectic form and are in general distinct from the de Rham cohomology, so that they provide new symplectic invariants. As shown in \cite{tsengyau} (cf. also Proposition~\ref{lefschetz} below), these new invariants actually agree with the de Rham cohomology if and only if the Hard Lefschetz property holds.

In this Section we discuss these cohomological invariants, proving that they can be computed using invariant forms, provided this is the case for the Rham cohomology (see Theorem~\ref{cinesi}).
This result will allow us to go through the list of symplectic structures on solvable Lie algebras (Table 3), to see which solvmanifolds satisfy the Hard Lefschetz property (Theorem~\ref{HL}).

\smallskip

Let $ (M,\omega) $ be a symplectic manifold of dimension $ 2n $,  one defines the Lefschetz operator 

\begin{displaymath}
\left. \begin{array}{l}
L:\Omega^{k}(M) \rightarrow \Omega^{k+2}(M) \\
\qquad \qquad \eta \mapsto \eta\wedge\omega
\end{array} \right.
\end{displaymath}
its dual operator $ \varLambda:\Omega^{k}(M) \rightarrow \Omega^{k-2}(M)$, and the symplectic star operator \linebreak $*_s: \Omega ^k(M)\rightarrow \Omega ^{2n-k}(M)$, such that for any $\gamma ,\beta \in  \Omega ^k(M)$, 
$$\gamma \wedge *_s\beta :=(\omega ^{-1})^k(\gamma ,\beta )dvol \, .
$$
\begin{remark} (see \cite{tsengyau})
\begin{enumerate}
\item $ *_s*_s=1 $.
\item Using coordinates $ (x_{1},..,x_{2n})$ on $M$ the above operators are defined in the following way: $$ \varLambda(\eta):=\frac{1}{2}(\omega^{-1})^{ij}i_{\partial_{x_{i}}}i_{\partial_{x_{j}}}\eta$$
where $ i $ is the interior product,
$$ \gamma \wedge *_s\beta := \frac{1}{k!}(\omega ^{-1})^{i_1j_1}(\omega ^{-1})^{i_2j_2} ... (\omega ^{-1})^{i_kj_k}\gamma _{i_1i_2...i_k}\beta _{j_1j_2...j_k}\frac{\omega ^n}{n!}.$$
\end{enumerate}
\end{remark}

\medskip

Using $ \varLambda $ one can construct two other differential operators: 
\begin{center}
 $d^{\wedge}:=(-1)^{k+1}*_sd*_s= d\varLambda-\varLambda d$ \; \; and\; \;  $ dd^{\wedge} $
\end{center}
with which one can define the following cohomologies 
\begin{displaymath}
H^{k}_{d^{\varLambda}}(M):= \dfrac{\ker d^{\varLambda}\cap\Omega^{k}(M)}{\mbox{im} \, d^{\varLambda}\cap\Omega^{k}(M)}
\end{displaymath}
 \begin{displaymath}
H^{k}_{d+d^{\varLambda}}(M):=\dfrac{\ker(d+d^{\varLambda})\cap\Omega^{k}(M)}{\mbox{im} \, dd^{\varLambda}\cap \Omega^{k}(M) } 
\end{displaymath}
\begin{displaymath}
H^{k}_{dd^{\varLambda}}(M):= \dfrac{\ker dd^{\varLambda}\cap \Omega^{k}(M)}{\mbox{im} \, d\cap\Omega^{k}(M) + \mbox{im} \,d^{\varLambda}\cap\Omega^{k}(M) }
\end{displaymath}
 \begin{displaymath}
 H^{k}_{d\cap d^{\varLambda}}(M):= H^{k}_{d}\cap H^{k}_{d^{\varLambda}}=H^{k}_{d+d^{\varLambda}}\cap H^{k}_{dd^{\varLambda}}=\dfrac{\ker(d+d^{\varLambda})\cap\Omega^{k}(M)}{\mbox{im} \, d\cap\Omega^{k}_{0}(M) + \mbox{im} \,d^{\varLambda}\cap \Omega^{k}_{0}(M) }
\end{displaymath}
where $ \Omega^{k}_{0}(M) $ is $ \ker dd^{\varLambda} \cap \Omega^{k}(M)$.

\smallskip

 We refer to \cite{tsengyau} for details and for the following propositions

\begin{proposition}\label{hodge}
\textnormal{(Tseng-Yau)} The operator $ *_{s} $ gives an isomorphism between $ H^{k}_{d}(M) $ and $ H^{2n-k}_{d^{\varLambda}} $ and between $ H^{k}_{d+d^{\varLambda}}(M) $ and $ H^{2n-k}_{dd^{\varLambda}} $.
\end{proposition} 

\begin{proposition}\label{lefschetz}
\textnormal{(Tseng-Yau)} On a compact symplectic manifold $ (M,\omega) $ the following properties are equivalent:
\begin{itemize}
\item the Hard Lefschetz property holds.
\item the canonical homomorphism $ H^{k}_{d+d^{\varLambda}}(M)\rightarrow H^{k}_{d}(M) $ is an isomorphism for all k.
\item the canonical homomorphism $ H^{k}_{d\cap d^{\varLambda}}(M)\rightarrow H^{k}_{d+d^{\varLambda}}(M) $ is an isomorphism for all k.
\end{itemize}
\end{proposition}

We are interested in the Lie groups associated to the six dimensional unimodular solvable non-nilpotent Lie algebras which admit a lattice and for which the de Rham cohomology can be computed by invariant forms, e.g., they are completely solvable. Indeed the following theorem holds:

\begin{theorem}\label{cinesi}
Let $ G $ be a Lie group admitting a left-invariant symplectic structure and a lattice $ \Gamma $ such that the quotient $ Q=G/\Gamma $ is compact. Let  $ \mathfrak{g} $ be the Lie algebra of $G$. 

If the inclusion
$\bigwedge^{*}(\mathfrak{g}^{*})\stackrel{i}{\hookrightarrow} \Omega^{*}(Q)$ is a quasi-isomorphism, (i.e., $ H^{*}_{d}(Q) \cong H^{*}_{d}(\mathfrak{g}) $), then $$ H^{*}_{d^{\varLambda}}(Q) \cong H^{*}_{d^{\varLambda}}(\mathfrak{g}), \quad H^{*}_{d+d^{\varLambda}}(Q) \cong H^{*}_{d+d^{\varLambda}}(\mathfrak{g}),$$ $$ H^{*}_{dd^{\varLambda}}(Q) \cong H^{*}_{dd^{\varLambda}}(\mathfrak{g}), \quad H^{*}_{d\cap d^{\varLambda}}(Q) \cong H^{*}_{d\cap d^{\varLambda}}(\mathfrak{g}).$$
\end{theorem}

\medskip

\begin{proof}
We divide the proof into four steps:

\noindent 1) \textit{We prove that the invariant cohomologies are well defined, i.e.,  the algebra of invariant forms $ \bigwedge^{*}(\mathfrak{g}^{*}) $ is closed for the operator $ d^{\varLambda} $.}
  
To this aim it suffices to prove that the operator $ *_s $ sends invariant forms to invariant forms. If $ \mathcal{L}:G\rightarrow G $ denotes the left translation, then $ \alpha $ and $ \beta $ are invariant if $ \mathcal{L}^* \alpha =\alpha $ and $ \mathcal{L}^* \beta =\beta $. Then
\[
\begin{array}{ll}
  &\mathcal{L}^*(\alpha \wedge *_s\beta ) = \mathcal{L}^*\left (\dfrac{1}{k!}(\omega ^{-1})^{i_1j_1}(\omega ^{-1})^{i_2j_2} ... (\omega ^{-1})^{i_kj_k}\alpha _{i_1i_2...i_k}\beta _{j_1j_2...j_k}\dfrac{\omega ^n}{n!}\right ) \\
  &{}\\
 &= \dfrac{1}{k!}(\mathcal{L}^*(\omega) ^{-1})^{i_1j_1}(\mathcal{L}^*(\omega) ^{-1})^{i_2j_2} ... (\mathcal{L}^*(\omega) ^{-1})^{i_kj_k}\mathcal{L}^*(\alpha _{i_1i_2...i_k})\mathcal{L}^*(\beta _{j_1j_2...j_k})\dfrac{\mathcal{L}^*(\omega ^n)}{n!} \\
 &{}\\
 &= \dfrac{1}{k!}(\omega ^{-1})^{i_1j_1}(\omega ^{-1})^{i_2j_2} ... (\omega ^{-1})^{i_kj_k}\alpha _{i_1i_2...i_k}\beta _{j_1j_2...j_k}\dfrac{\omega ^n}{n!}=\alpha \wedge *_s\beta \, .
\end{array}
\]

Therefore, $\alpha \wedge *_s\beta =\mathcal{L}^*(\alpha \wedge *_s\beta) =\mathcal{L}^*(\alpha )\wedge \mathcal{L}^*(*_s\beta )=\alpha \wedge \mathcal{L}^*(*_s\beta )$ and so $$*_s\beta =\mathcal{L}^*(*_s\beta )\, .
$$
\noindent 2) \textit{We show that $H^{*}_{d^{\varLambda}}(Q) \cong H^{*}_{d^{\varLambda}}(\mathfrak{g})$,  $H^{*}_{d\cap d^{\varLambda}}(Q) \cong H^{*}_{d\cap d^{\varLambda}}(\mathfrak{g})$
and that $H^{*}_{d+d^{\varLambda}}(Q) \cong H^{*}_{d+d^{\varLambda}}(\mathfrak{g})$ if and only if $H^{*}_{dd^{\varLambda}}(Q) \cong H^{*}_{dd^{\varLambda}}(\mathfrak{g})$. }

\smallskip

Notice that 1) and Proposition~\ref{hodge} imply the commutativity of the diagram
\begin{displaymath}
\xymatrix{
H^k_d(\mathfrak{g}) \ar[r]^{*_s}_\sim   \ar[d]^\sim   &
H^{2n-k}_{d^\wedge }(\mathfrak{g}) \ar[d]   \\
H^k_d(Q) \ar[r]_{*_s}^\sim  & H^{2n-k}_{d^\wedge } (Q)}
\end{displaymath}
so that, since by assumption the isomorphism  holds for $H_{d}$, it holds for $H_{d^\varLambda}$, i.e., $H^{*}_{d^{\varLambda}}(Q) \cong H^{*}_{d^{\varLambda}}(\mathfrak{g})$.
Moreover, since $H^{k}_{d\cap d^{\varLambda}}(Q):= H^{k}_{d}\cap H^{k}_{d^{\varLambda}}$,  the isomorphism holds also for the $d\cap d^{\varLambda}$-cohomology. 

Hence Proposition~\ref{hodge} implies that if the isomorphism between cohomology and invariant cohomology holds for $H_{d+d^{\varLambda}}$, then it is also true  for $ H_{dd^{\varLambda}} $ and vice versa.\\
\medskip

\noindent 3)  \textit{$ i^{*}: H^*_{d+d^{\varLambda}}(\mathfrak{g}^{*})\rightarrow H^*_{d+d^{\varLambda}}(Q) $ is injective}, (cf. \cite[page 123]{rag}). 

Since $Q$ is compact, there exists an invariant metric $\langle \, , \, \rangle$ on $Q$. One can use this metric to define the adjoint operators of $d, \, d^{\varLambda}, \, d+d^{\varLambda}$ and $dd^{\varLambda}$. Let $ \bigwedge^{\bot k}(\mathfrak{g}^{*}) $ be the orthogonal complement of $ \bigwedge^{k}(\mathfrak{g}^{*}) $ in $ \Omega^{k}(Q) $. Then $ \Omega^{k}(Q)=\bigwedge^{k}(\mathfrak{g}^{*})\oplus \bigwedge^{\bot k}(\mathfrak{g}^{*}) $ and $ \bigwedge^{k}(\mathfrak{g}^{*}) $ and $ \bigwedge^{\bot k}(\mathfrak{g}^{*})  $ are closed under $ d+d^{\varLambda} $ and $ dd^{\varLambda} $.

If $i^*[\alpha ]:=[i(\alpha )]=0$, then there exists a form $\eta \in \Omega(Q) $ such that 
$$i(\alpha )=dd^\wedge \eta =dd^\wedge (\tilde \eta + \tilde \eta ^\bot )=dd^\wedge \tilde \eta + dd^\wedge \tilde \eta ^\bot, $$ with $\tilde \eta\in \bigwedge^{k}(\mathfrak{g}^{*}) $ and  $\tilde \eta ^\bot \in \bigwedge^{\bot k}(\mathfrak{g}^{*})$. 

Moreover $dd^\wedge \tilde \eta \in \bigwedge^{k}(\mathfrak{g}^{*}) $, so $i(\alpha - dd^\wedge \tilde \eta)= dd^\wedge \tilde \eta ^\bot$ and $[\alpha ]=[\alpha - dd^\wedge \tilde \eta]$. 

So we can choose $\tilde \alpha :=\alpha - dd^\wedge \tilde \eta$ as a representative of the cohomology class $[\alpha]$ in $H^*_{d+d^{\varLambda}}(\mathfrak{g}^{*})$. 

Observe that $\tilde \alpha \in \bigwedge^{k}(\mathfrak{g}^{*}) $ so $(dd^\wedge )^*\tilde \alpha \in \bigwedge^{k}(\mathfrak{g}^{*})$ and $i((dd^\wedge )^*\tilde \alpha)=(dd^\wedge )^*i(\tilde \alpha )=(dd^\wedge )^*dd^\wedge \tilde \eta ^\bot\in \bigwedge^{k}(\mathfrak{g}^{*}) $, but then $\tilde \eta ^\bot \in \bigwedge^{\bot k}(\mathfrak{g}^{*})$ is orthogonal to $(dd^\wedge )^*dd^\wedge \tilde \eta ^\bot \in  \bigwedge^{k}(\mathfrak{g}^{*})$. This implies 
$$0=\langle \tilde \eta ^\bot ,(dd^\wedge )^*dd^\wedge \tilde \eta ^\bot \rangle=\langle dd^\wedge \tilde \eta ^\bot ,dd^\wedge \tilde \eta ^\bot \rangle\, ,
$$ so $dd^\wedge \tilde \eta ^\bot =0$. But then $i(\alpha )=dd^\wedge \tilde \eta $, with $\tilde \eta $ in $  \bigwedge^{k}(\mathfrak{g}^{*}) $, so $\alpha =dd^\wedge \tilde \eta $ in $ \bigwedge^{k}(\mathfrak{g}^{*}) $, that is $[\alpha ]=0$ belongs to $ H^*_{d+d^\wedge }(\mathfrak{g}^{*})$.

\begin{remark} We can similarly prove that also $ i^{*}: H^*_{dd^{\varLambda}}(\mathfrak{g}^{*})\rightarrow H^*_{dd^{\varLambda}}(Q) $ is injective.\\In particular 3) is always true, independent on the fact that the map $i$ is a quasi-isomorphism.
\end{remark}

\noindent 4)  
\textit{$i^{*}: H^*_{d+d^{\varLambda}}(\mathfrak{g}^{*})\rightarrow H^*_{d+d^{\varLambda}}(Q)$ is surjective}. 

Let $ \eta \in \Omega^{k}(Q) $ be such that $ d\eta=d^{\varLambda}\eta=0 $. Then the cohomology class $ [\eta]^{k}_{d+d^{\varLambda}} $ is well defined. But also $ [\eta]^{k}_{d} $ and $  [\eta]^{k}_{d^{\varLambda}} $ exist and by hypothesis they have an invariant representative: $ \eta=\tilde{\eta}_{1} +d\mu_{1} $ and $ \eta=\tilde{\eta}_{2} +d^{\varLambda}\mu_{1} $ with $ \eta_{1},\eta_{2} \in \bigwedge^{*}(\mathfrak{g}^{*})$ and $ d\tilde{\eta}_{1}=d^{\varLambda}\tilde{\eta}_{2}=0 $.

Since $ dd^{\varLambda}\eta=0 $, the cohomology class $ [\eta]^{k}_{dd^{\varLambda}} $ exists and 
$$ \eta=\frac{1}{2}(\tilde{\eta}_{1}+\tilde{\eta}_{2})+ d\frac{\mu_{1}}{2}+ d^{\varLambda}\frac{\mu_{2}}{2} \, $$
then $ \frac{1}{2}(\tilde{\eta}_{1}+\tilde{\eta}_{2}) $ is an invariant representative for $ [\eta]^{k}_{dd^{\varLambda}} $.

Now we apply the isomorphism of Proposition~\ref{hodge}: 
$$[*_s\eta]^{2n-k}_{d+d^{\varLambda}}\cong [\eta]^{k}_{dd^{\varLambda}}= \left[ \dfrac{\tilde{\eta}_{1}+\tilde{\eta}_{2}}{2}\right]^{k}_{dd^{\varLambda}}\cong \left[ *_s \big( \dfrac{\tilde{\eta}_{1}+\tilde{\eta}_{2}}{2}\big)\right]^{2n-k}_{d+d^{\varLambda}}. $$
Let $ *_s\eta=N, \; *_s\tilde{\eta}_{1}=N_{1}, \; *_s\tilde{\eta}_{2}=N_{2} $. Then $ \frac{N_1+N_2}{2} $ is an invariant representative in $ [N]^{2n-k}_{d+d^{\varLambda}} $.

To complete the proof we have to show that every $ N \in \Omega^{2n-k}(Q) $ such that $ dN=d^{\varLambda}N=0 $ is of the form $ N=*_s\eta $ with $ \eta \in \Omega^{k}(Q) $ and $ d\eta=d^{\varLambda}\eta=0$. 

To this aim, it is sufficient to impose $ \eta:=*_sN $, then $ *_s\eta=*_s*_sN=N $. Moreover $ d^{\wedge}:=(-1)^{k+1}*_sd*_s $, so $ *_sd^{\varLambda}=(-1)^{k+1}d*_s $ and $ d^{\varLambda}*_s=(-1)^{k+1}*_sd $. Then for every $ \beta \in \Omega^{k}(Q) $ if $ d^{\varLambda}\beta=0 $, also $ *_sd^{\varLambda}\beta=0 $ and then $ d*_s\beta=0 $ and similarly if $ d\beta=0 $, then $ d^{\varLambda}*_s\beta=0 $. 

Hence $ d\eta=d^{\varLambda}\eta=0 $.

\end{proof}

\begin{remark}
In particular Theorem~\ref{cinesi} applies in the following cases:
\begin{itemize}
\item If $G$ is nilpotent, using Nomizu theorem, \cite{nomizu}.
\item If $G$ is completely solvable, using Hattori theorem, \cite{hattori}.
\item If $ \mbox{Ad}(G) $ and $ \mbox{Ad}(\Gamma) $ have the same algebraic closure, using Mostow theorem, \cite{rag}, cf. also \cite{guan2, cf}.
\end{itemize}
\end{remark}

Using Theorem~\ref{cinesi}  and Proposition~\ref{lefschetz} we can examine which symplectic solvmanifold whose Lie algebra is in Table 2 with $G$ completely solvable, is Hard Lefschetz.

Let $ \{\alpha_{1}, ..., \alpha_{6}\}$  be the dual basis of $ \{X_1,...,X_6\} $. Then a generic element in $ \bigwedge \mathfrak{g}^{*2} $ is $ \beta = \sum_{i<j}b_{i,j}\alpha_{ij} $, where we use the notation $ \alpha_{i_{1}...i_{n}} := \alpha_{i_{1}} \wedge ... \wedge \alpha_{i_{n}} $. 

For any such solvmanifold we perform the computation only for a particular choice of the symplectic form. Namely we consider the form composed by the fewest possible generators $ \alpha_{ij} $ of $ \bigwedge^{2}(\mathfrak{g}^{*}) $, and we check if the Hard Lefschetz property holds only for this particular choice. This is because computations are very involved for a generic symplectic form.

\medskip

The symplectic and completely solvable Lie algebras in Table 2 whose Lie group admits a lattice are the following, \cite{bock}:
\medskip
\\ $ \mathfrak{g}_{3.1}\oplus 3\mathbb{R}, \quad \mathfrak{g}_{3.1}\oplus\mathfrak{g}_{3.4}, \quad \mathfrak{g}_{3.4}\oplus\mathfrak{g}_{3.4}, \quad \mathfrak{g}_{5.7}^{p,-p,-1}\oplus \mathbb{R}, \quad \mathfrak{g}_{5.8}\oplus \mathbb{R}, \quad  \mathfrak{g}_{5.15}\oplus \mathbb{R}, \quad \mathfrak{g}_{6.3}, \newline \mathfrak{g}_{6.15}, \quad \mathfrak{g}_{6.21}^{0}, \quad \mathfrak{g}_{6.23}^{0,0,\pm 1}, \quad \mathfrak{g}_{6.29}^{0,0,\pm 1}, \quad \mathfrak{g}_{6.29}^{0,0,0}, \quad \mathfrak{g}_{6.54}^{0,-1}, \quad \mathfrak{g}_{6.78} $.
\bigskip

Next we list only the cases when Hard Lefschetz property holds.

\begin{itemize}
\item $ \mathfrak{g}_{3.4}\oplus 3\mathbb{R} : \; $
\begin{tabular}{@{} l @{} }
$\omega=\omega_{1,2}\alpha_{12}+\omega_{3,6}\alpha_{36}+\omega_{4,5}\alpha_{45}, \; \tilde  \omega=\omega_{1,2}\alpha_{12}+\omega_{3,4}\alpha_{34}+\omega_{5,6}\alpha_{56}$ \\
$\hat \omega=\omega_{1,2}\alpha_{12}+\omega_{3,5}\alpha_{35}+\omega_{4,4}\alpha_{46}$ \\
\end{tabular}
 \medskip \\
\begin{tabular}{ @{} c @{} c }
$b^{1}_{d}= b^{1}_{d+d^{\varLambda}}= b^{1}_{d\cap d^{\varLambda}}=4$  \\
$b^{2}_{d}= b^{2}_{d+d^{\varLambda}}= b^{2}_{d\cap d^{\varLambda}}=7$  \\
$ b^{3}_{d}= b^{3}_{d+d^{\varLambda}}= b^{3}_{d\cap d^{\varLambda}}=8 $   \\
\end{tabular}
\bigskip
\item $ \mathfrak{g}_{3.4}\oplus\mathfrak{g}_{3.4} : \quad \omega=\omega_{1,2}\alpha_{12}+\omega_{3,6}\alpha_{36}+\omega_{4,5}\alpha_{45}$ \medskip \\
\begin{tabular}{ @{} c @{} c  }
$b^{1}_{d}= b^{1}_{d+d^{\varLambda}}= b^{1}_{d\cap d^{\varLambda}}=2$ \\
$b^{2}_{d}= b^{2}_{d+d^{\varLambda}}= b^{2}_{d\cap d^{\varLambda}}=3$   \\
$ b^{3}_{d}= b^{3}_{d+d^{\varLambda}}= b^{3}_{d\cap d^{\varLambda}}=4 $   
\end{tabular}
\bigskip
\item $ \mathfrak{g}_{5.7}^{p,-p,-1}\oplus \mathbb{R} :  $
\begin{itemize}
\item[$p=1$:] $ \; \omega=\omega_{1,4}\alpha_{14}+\omega_{2,3}\alpha_{23}+\omega_{5,6}\alpha_{56}, \quad \tilde \omega=\omega_{1,3}\alpha_{13}+\omega_{2,4}\alpha_{24}+\omega_{5,6}\alpha_{56} $
\begin{tabular}{ @{} c @{} c }
$b^{1}_{d}= b^{1}_{d+d^{\varLambda}}= b^{1}_{d\cap d^{\varLambda}}=2$  \\
$b^{2}_{d}= b^{2}_{d+d^{\varLambda}}= b^{2}_{d\cap d^{\varLambda}}=5$  \\
$ b^{3}_{d}= b^{3}_{d+d^{\varLambda}}= b^{3}_{d\cap d^{\varLambda}}=8 $
\end{tabular} \medskip
\item[$p\neq 1$:] $\; \omega=\omega_{1,4}\alpha_{14}+\omega_{2,3}\alpha_{23}+\omega_{5,6}\alpha_{56}  $ \\
\begin{tabular}{ @{} c @{} c }
$b^{1}_{d}= b^{1}_{d+d^{\varLambda}}= b^{1}_{d\cap d^{\varLambda}}=2$  \\
$b^{2}_{d}= b^{2}_{d+d^{\varLambda}}= b^{2}_{d\cap d^{\varLambda}}=3$   \\
$ b^{3}_{d}= b^{3}_{d+d^{\varLambda}}= b^{3}_{d\cap d^{\varLambda}}=4 $   
\end{tabular}
\end{itemize}
\end{itemize}

\bigskip 
We have then proved:

\begin{theorem}\label{HL}
There exists a symplectic structure for which the following solvmanifolds are Hard Lefschetz:
$$ (G_{5.7}^{p,-p,-1}\otimes\mathbb{R})/\Gamma, \quad (G_{3.4}\otimes 3\mathbb{R})/\Gamma, \quad (G_{3.4}\otimes G_{3.4})/\Gamma\, .$$
\end{theorem}

\begin{remark} The case of $  (G_{5.7}^{p,-p,-1}\otimes\mathbb{R})/\Gamma $ was already considered in \cite{bock}.
\end{remark}
\newpage

\section{Appendix}
\section*{Table 2:  {\rm Six dimensional solvable (non nilpotent) unimodular Lie algebras}}
\begin{center}
{\small
\begin{tabular}{|c|c|}
\hline
Algebra & Structure equations \\
\hline
$ \mathfrak{g}_{6.1}^{a,b,c,e}$ & $ d\alpha_1=-\alpha_{16}, \; d\alpha_2=-a\alpha_{26}, \; d\alpha_3=-b\alpha_{36} $ \\ { \scriptsize$0< |e|\leq |c| \leq |b| \leq |a| \leq 1$,}  & $ d\alpha_4=-c\alpha_{46},  \; d\alpha_5=-e\alpha_{56}, \; d\alpha_6 =0 $ \\ {\scriptsize $a+b+c+e=-1$} &  \\
\hline 
$ \mathfrak{g}_{6.2}^{a,c,e}$ & $ d\alpha_1=-a\alpha_{16}-\alpha_{26},\; d\alpha_2=-a\alpha_{26},\; d\alpha_3=-\alpha_{36}, $ \\ {\scriptsize $0< |e|\leq |c| \leq 1$, } & $ d\alpha_4=-c\alpha_{46},  \; d\alpha_5=-e\alpha_{56}, \; d\alpha_6= 0$ \\ {\scriptsize$2a+c+e=-1$ }&  \\
\hline
$ \mathfrak{g}_{6.3}^{-\frac{a+1}{3},a}$ & $ d\alpha_1=\frac{a+1}{3}\alpha_{16} - \alpha_{26},\; d\alpha_2=\frac{a+1}{3}\alpha_{26}-\alpha_{36},  $ \\ {\scriptsize $0< |a|\leq 1$}  & $ d\alpha_3=\frac{a+1}{3}\alpha_{36},\;d\alpha_4=-\alpha_{46}, $ \\ & $ d\alpha_5=-a\alpha_{56},\; d\alpha_6=0 $ \\
\hline 
$ \mathfrak{g}_{6.4}^{-\frac{1}{4}}$ & $ d\alpha_1=\frac{1}{4}\alpha_{16} - \alpha_{26},\; d\alpha_2=  \frac{1}{4}\alpha_{26}-\alpha_{36}, $ \\  & $ d\alpha_3=\frac{1}{4}\alpha_{36}-\alpha_{46},\; d\alpha_4=\frac{1}{4}\alpha_{46}, \; d\alpha_5=-\alpha_{56},\; d\alpha_6= 0$\\
\hline 
$ \mathfrak{g}_{6.6}^{a,b}$ & $ d\alpha_1=-\alpha_{16},\; d\alpha_2=-a\alpha_{26}-\alpha_{36},\; d\alpha_3=-a\alpha_{36},$ \\  {\scriptsize$ a\leq b, \quad a+b=-\frac{1}{2}$ } & $  d\alpha_4=-b\alpha_{46}-\alpha_{56}, \;d\alpha_5= -b\alpha_{56},d\alpha_6 =0  $ \\
\hline
$ \mathfrak{g}_{6.7}^{a,-\frac{2}{3}a}$ & $ d\alpha_1=-a\alpha_{16} - \alpha_{26},\; d\alpha_2=-a\alpha_{26}-\alpha_{36}, \; d\alpha_3=-a\alpha_{36}, $ \\ {\scriptsize$ a\neq 0 $} & $ d\alpha_4= \frac{2}{3}a\alpha_{46}-\alpha_{56}, \; d\alpha_5= \frac{2}{3}a\alpha_{56}, d\alpha_6=0 $ \\
\hline 
$ \mathfrak{g}_{6.8}^{a,b,c,p}$ & $ d\alpha_1=-a\alpha_{16},\; d\alpha_2=-b\alpha_{26},\; d\alpha_3=-c\alpha_{36}, $\\ {\scriptsize $0< |c| \leq |b| \leq |a|, $ } & $ d\alpha_4= -p\alpha_{46}-\alpha_{56},\; d\alpha_5=\alpha_{46} -p\alpha_{56}, \; d\alpha_6= 0 $  \\ {\scriptsize $ a+b+c+2p=0 $} & \\
\hline 
$ \mathfrak{g}_{6.9}^{a,b,p}$ & $ d\alpha_1=-a\alpha_{16}, \; d\alpha_2=-b\alpha_{26}-\alpha_{36}, \; d\alpha_3=-b\alpha_{36}, $\\  {\scriptsize$a\neq 0, \; a+2b+2p=0$}  & $ d\alpha_4= -p\alpha_{46}-\alpha_{56}, d\alpha_5=\alpha_{46} -p\alpha_{56}, d\alpha_6=0 $  \\
\hline 
$ \mathfrak{g}_{6.10}^{a,-\frac{3}{2}a}$ & $ d\alpha_1=-a\alpha_{16} - \alpha_{26},\; d\alpha_2=-a\alpha_{26}-\alpha_{36},\; d\alpha_3=-a\alpha_{36},  $ \\ &$d\alpha_4=\frac{3}{2}a\alpha_{46}-\alpha_{56}, \; d\alpha_5=\alpha_{46}+\frac{3}{2}a\alpha_{56},d\alpha_6= 0$ \\
\hline 
$ \mathfrak{g}_{6.11}^{a,p,q,s}$ & $ d\alpha_1=-a\alpha_{16}, \; d\alpha_2=-p\alpha_{26}-\alpha_{36},\; d\alpha_3=\alpha_{26}-p\alpha_{36}, $ \\ {\scriptsize$ as\neq 0, \; a+2p+2q=0 $} & $  d\alpha_4=-q\alpha_{46}-s\alpha_{56},\; d\alpha_5=s\alpha_{46}-q\alpha_{56},\; d\alpha_6= 0 $\\
\hline 
$ \mathfrak{g}_{6.12}^{-4p,p}$ & $ d\alpha_1=4p\alpha_{16},\; d\alpha_2=-p\alpha_{26}-\alpha_{36}-\alpha_{46},   $ \\ {\scriptsize$ p\neq 0$} & $ d\alpha_3=\alpha_{26}-p\alpha_{36}-\alpha_{56},\;d\alpha_4=-p\alpha_{46}-\alpha_{56}, $\\ & $ d\alpha_5=\alpha_{46}-p\alpha_{56},\; d\alpha_6= 0 $ \\
\hline 
$ \mathfrak{g}_{6.13}^{a,b,h}$ & $ d\alpha_1=-\alpha_{23}-(a+b)\alpha_{16},\; d\alpha_2=-a\alpha_{26},\; d\alpha_3=-b\alpha_{36},  $ \\{\scriptsize $ a\neq 0, \; 2a+2b+h=-1 $} & $ d\alpha_4=-\alpha_{46},\; d\alpha_5=-h\alpha_{56},\; d\alpha_6= 0 $\\
\hline 
$ \mathfrak{g}_{6.14}^{a,b}$ & $ d\alpha_1=-\alpha_{23}-(a+b)\alpha_{16}-\alpha_{56},\; d\alpha_2=-a\alpha_{26},$ \\ {\scriptsize$ a\neq 0, \quad a+b=-\frac{1}{3} $ } & $ d\alpha_3=-b\alpha_{36},\; d\alpha_4=-\alpha_{46},\; d\alpha_5=-(a+b)\alpha_{56}, \; d\alpha_6= 0  $\\
\hline 
$ \mathfrak{g}_{6.15}^{-1}$ & $ d\alpha_1=-\alpha_{23},\; d\alpha_2=-\alpha_{26},\; d\alpha_3=\alpha_{36}, $ \\ & $ d\alpha_4=-\alpha_{26}-\alpha_{46},\; d\alpha_5=-\alpha_{36}+\alpha_{56},\; d\alpha_6=0 $\\
\hline 
$ \mathfrak{g}_{6.17}^{-\frac{1}{2},0}$ & $ d\alpha_1=-\alpha_{23}+\frac{1}{2}\alpha_{16},\; d\alpha_2=\frac{1}{2}\alpha_{26}, \; d\alpha_3= 0, $ \\ & $ d\alpha_4=-\alpha_{36},\; d\alpha_5=-\alpha_{56},\; d\alpha_6= 0 $ \\
\hline 
$ \mathfrak{g}_{6.18}^{a,-2a-3}$ & $ d\alpha_1=-\alpha_{23}-(1+a)\alpha_{16},\; d\alpha_2=-a\alpha_{26},\; d\alpha_3=-\alpha_{36},$ \\ {\scriptsize$ a \neq -\frac{3}{2} $ }& $ d\alpha_4=-\alpha_{36}-\alpha_{46},\; d\alpha_5=(2a+3)\alpha_{56},\; d\alpha_6= 0 $\\
\hline 
$ \mathfrak{g}_{6.19}^{-\frac{4}{3}}$ & $ d\alpha_1=-\alpha_{23}+\frac{1}{3}\alpha_{16}-\alpha_{56},\; d\alpha_2=\frac{4}{3}\alpha_{26},\; d\alpha_3=-\alpha_{36}, $ \\ & $d\alpha_4= -\alpha_{36}-\alpha_{46},\; d\alpha_5=\frac{1}{3}\alpha_{56},\; d\alpha_6= 0 $\\ 
\hline 
$ \mathfrak{g}_{6.20}^{-3}$ & $ d\alpha_1=-\alpha_{23}-\alpha_{16}-\alpha_{46}, \; d\alpha_2=0,\; d\alpha_3= -\alpha_{36}, $ \\ & $ d\alpha_4=-\alpha_{36}-\alpha_{46},\; d\alpha_5=3\alpha_{56},\; d\alpha_6= 0 $ \\ 
\hline 
$ \mathfrak{g}_{6.21}^{a}$ & $ d\alpha_1=-\alpha_{23}-2a\alpha_{16},\; d\alpha_2=-a\alpha_{26},\; d\alpha_3=-\alpha_{26}-a\alpha_{36}, $ \\{\scriptsize $ a \neq -\frac{1}{4} $} & $ d\alpha_4=-\alpha_{46},\; d\alpha_5=(4a+1)\alpha_{56},\; d\alpha_6=  0$\\
\hline 
$ \mathfrak{g}_{6.22}^{-\frac{1}{6}}$ & $ d\alpha_1=-\alpha_{23}+\frac{1}{3}\alpha_{16}-\alpha_{56},\; d\alpha_2=\frac{1}{6}\alpha_{26}, $  \\ & $ d\alpha_3=-\alpha_{26}+\frac{1}{6}\alpha_{36},\; d\alpha_4=-\alpha_{46}, $ \\ & $  d\alpha_5=\frac{1}{3}\alpha_{56}, \; d\alpha_6=0 $ \\
\hline 

\end{tabular}
}
\end{center}

\begin{center}
{\small
\begin{tabular}{|c|c|}
\hline
$ \mathfrak{g}_{6.23}^{a,-7a,\varepsilon}$ & $ d\alpha_1=-\alpha_{23}-2a\alpha_{16}-\varepsilon\alpha_{56},\; d\alpha_2=-a\alpha_{26},  $ \\{\scriptsize $ a\varepsilon =0 $} & $d\alpha_3=-\alpha_{26}-a\alpha_{36},\; d\alpha_4=-\alpha_{36}-a\alpha_{46},$ \\ & $  d\alpha_5=5a\alpha_{56}, \; d\alpha_6=0 $ \\
\hline 
$ \mathfrak{g}_{6.25}^{b,-1-b}$ & $ d\alpha_1=-\alpha_{23}+b\alpha_{16},\; d\alpha_2=-\alpha_{26},\; d\alpha_3=(1+b)\alpha_{36}, $ \\  & $ d\alpha_4=-b\alpha_{46},\;d\alpha_5= -\alpha_{46}-b\alpha_{56},\;  d\alpha_{6}=0$ \\

\hline
$ \mathfrak{g}_{6.26}^{-1}$ & $ d\alpha_1=-\alpha_{23}-\alpha_{56},\; d\alpha_2=-\alpha_{26},\; d\alpha_3=\alpha_{36}, $ \\  & $ d\alpha_4=0, \; d\alpha_5=-\alpha_{46},\; d\alpha_6= 0 $\\
\hline
$ \mathfrak{g}_{6.27}^{-2b,b}$ & $ d\alpha_1=-\alpha_{23}+b\alpha_{16},\;  d\alpha_2=2b\alpha_{26},\;d\alpha_3=-b\alpha_{36}, $ \\{\scriptsize $ b\neq 0 $} & $ d\alpha_4=-\alpha_{36}-b\alpha_{46},\; d\alpha_5=-\alpha_{46}-b\alpha_{56},\; d\alpha_6= 0$ \\
\hline
$ \mathfrak{g}_{6.28}^{-2}$ & $ d\alpha_1=-\alpha_{23}-2\alpha_{16},\; d\alpha_2=-\alpha_{26},\; d\alpha_3=-\alpha_{26}-\alpha_{36}, $ \\ & $d\alpha_4= 2\alpha_{46},\; d\alpha_5=-\alpha_{46}+2\alpha_{56},\;d\alpha_6= 0$ \\
\hline
$ \mathfrak{g}_{6.29}^{-2b,b,\varepsilon}$ & $ d\alpha_1=-\alpha_{23}+b\alpha_{16}-\varepsilon\alpha_{56},\; d\alpha_2=2b\alpha_{26},\;d\alpha_3=-b\alpha_{36}, $ \\{\scriptsize $ b\varepsilon =0$ }& $ d\alpha_4=-\alpha_{36}-b\alpha_{46},\; d\alpha_5=-\alpha_{46}-b\alpha_{56},\; d\alpha_6= 0 $\\
\hline
$ \mathfrak{g}_{6.32}^{a,-6a-h,h,\varepsilon}$ & $ d\alpha_1=-\alpha_{23}-2a\alpha_{16}-\varepsilon\alpha_{46}, \; d\alpha_2= -a\alpha_{26}+\alpha_{36},   $ \\ {\scriptsize$ a>-\frac{1}{4}h, \varepsilon h =0 $}  &  $d\alpha_3=-\alpha_{26}-a\alpha_{36},\; d\alpha_4=-(2a+h)\alpha_{46}, $ \\ & $  d\alpha_5=(6a+h)\alpha_{56}, \; d\alpha_6 =0 $ \\
\hline
$ \mathfrak{g}_{6.33}^{a,-6a}$ & $ d\alpha_1=-\alpha_{23}-2a\alpha_{16}-\alpha_{56}, \; d\alpha_2= -a\alpha_{26}+\alpha_{36},  $ \\ {\scriptsize$ a\geq 0 $}  &  $ d\alpha_3=-\alpha_{26}-a\alpha_{36},\; d\alpha_4=6a\alpha_{46}, $ \\ & $  d\alpha_5=-2a\alpha_{56}, \; d\alpha_6 =0 $ \\
\hline
$ \mathfrak{g}_{6.34}^{a,-4a,\varepsilon}$ & $ d\alpha_1=-\alpha_{23}-2a\alpha_{16}-\varepsilon\alpha_{56}, \; d\alpha_2= -a\alpha_{26}+\alpha_{36},   $ \\ {\scriptsize$ \varepsilon a=0 $}  &  $  d\alpha_3=-\alpha_{26}-a\alpha_{36},\;d\alpha_4=2a\alpha_{46}, \; d\alpha_5=2a\alpha_{56}, \; d\alpha_6 =0 $ \\
\hline
$ \mathfrak{g}_{6.35}^{a,b,c}$ & $ d\alpha_1=-\alpha_{23}-(a+b)\alpha_{16}, \; d\alpha_2= -a\alpha_{26},\; d\alpha_3=-b\alpha_{36},   $ \\ {\scriptsize$ a+b+c=0, \; a^2+b^2\neq 0 $}  &  $ d\alpha_4=-c\alpha_{46}+\alpha_{56}, \; d\alpha_5=-\alpha_{46}-c\alpha_{56}, \; d\alpha_6 =0 $ \\
\hline
$ \mathfrak{g}_{6.36}^{a,-2a}$ & $ d\alpha_1=-\alpha_{23}-2a\alpha_{16}, \; d\alpha_2= -a\alpha_{26},\; d\alpha_3=-\alpha_{26}-a\alpha_{36},   $ \\   &  $ d\alpha_4=2a\alpha_{46}+\alpha_{56}, \; d\alpha_5=-\alpha_{46}+2a\alpha_{56}, \; d\alpha_6 =0 $ \\
\hline
$ \mathfrak{g}_{6.37}^{a,-2a,s}$ & $ d\alpha_1=-\alpha_{23}-2a\alpha_{16}, \; d\alpha_2= -a\alpha_{26}+\alpha_{36},\; d\alpha_3=-\alpha_{26}-a\alpha_{36},   $ \\ {\scriptsize$ s\neq 0 $}  &  $ d\alpha_4=2a\alpha_{46}+s\alpha_{56}, \; d\alpha_5=-s\alpha_{46}+2a\alpha_{56}, \; d\alpha_6 =0 $ \\
\hline
$ \mathfrak{g}_{6.38}^{0}$ & $ d\alpha_1=-\alpha_{23}, \; d\alpha_2=\alpha_{36},\; d\alpha_3=-\alpha_{26},   $ \\   &  $ d\alpha_4=-\alpha_{26}+\alpha_{56}, \; d\alpha_5=-\alpha_{36}-\alpha_{46}, \; d\alpha_6 =0 $ \\
\hline
$ \mathfrak{g}_{6.39}^{-4-3h,h}$ & $ d\alpha_1=-\alpha_{45}-(1+h)\alpha_{16}, \; d\alpha_2= -\alpha_{15}-(2+h)\alpha_{26},  $ \\ {\scriptsize$ h\neq -\frac{4}{3} $}  &  $ d\alpha_3=(4+3h)\alpha_{36}, \; d\alpha_4=-h\alpha_{46}, \; d\alpha_5=-\alpha_{56}, \; d\alpha_6 =0 $ \\
\hline
$ \mathfrak{g}_{6.40}^{-\frac{3}{2}}$ & $ d\alpha_1=-\alpha_{45}+\frac{1}{2}\alpha_{16}, \; d\alpha_2= -\alpha_{15}-\frac{1}{2}\alpha_{26}-\alpha_{36},  $ \\   &  $ d\alpha_3=-\frac{1}{2}\alpha_{36}, \; d\alpha_4=\frac{3}{2}\alpha_{46}, \; d\alpha_5=-\alpha_{56}, \; d\alpha_6 =0 $ \\
\hline
$ \mathfrak{g}_{6.41}^{-1}$ & $ d\alpha_1=-\alpha_{45}, \; d\alpha_2= -\alpha_{15}-\alpha_{26},\;d\alpha_3=\alpha_{36}-\alpha_{46},  $ \\   &  $ d\alpha_4=\alpha_{46}, \; d\alpha_5=-\alpha_{56}, \; d\alpha_6 =0 $ \\
\hline
$ \mathfrak{g}_{6.42}^{-\frac{5}{3}}$ & $ d\alpha_1=-\alpha_{45}+\frac{2}{3}\alpha_{16}, \; d\alpha_2= -\alpha_{15}-\frac{1}{3}\alpha_{26},\;d\alpha_3=\alpha_{36}-\alpha_{56},  $ \\   &  $ d\alpha_4=\frac{5}{3}\alpha_{46}, \; d\alpha_5=-\alpha_{56}, \; d\alpha_6 =0 $ \\
\hline
$ \mathfrak{g}_{6.44}^{-7}$ & $ d\alpha_1=-\alpha_{45}-2\alpha_{16}, \; d\alpha_2= -\alpha_{15}-3\alpha_{26},\;d\alpha_3=7\alpha_{36},  $ \\   &  $ d\alpha_4=-\alpha_{46}-\alpha_{56}, \; d\alpha_5=-\alpha_{56}, \; d\alpha_6 =0 $ \\
\hline
$ \mathfrak{g}_{6.47}^{-3,\varepsilon}$ & $ d\alpha_1=-\alpha_{45}-\alpha_{16}, \; d\alpha_2= -\alpha_{15}-\alpha_{26}-\varepsilon\alpha_{46},  $ \\{\scriptsize $ \varepsilon=0,\pm 1 $ } &  $d\alpha_3=3\alpha_{36},\; d\alpha_4=-\alpha_{46}, \; d\alpha_5=0, \; d\alpha_6 =0 $ \\
\hline
$ \mathfrak{g}_{6.54}^{2(1+l),l}$ & $ d\alpha_1=-\alpha_{35}-\alpha_{16}, \; d\alpha_2= -\alpha_{45}-l\alpha_{26},\;d\alpha_3=(1+2l)\alpha_{36},  $ \\ &  $ d\alpha_4=(2+l)\alpha_{46}, \; d\alpha_5=-2(1+l)\alpha_{56}, \; d\alpha_6 =0 $ \\
\hline
$ \mathfrak{g}_{6.55}^{-4}$ & $ d\alpha_1=-\alpha_{35}-\alpha_{16}-\alpha_{46}, \; d\alpha_2= -\alpha_{45}+3\alpha_{26}, $ \\ &  $d\alpha_3=-5\alpha_{36},\;  d\alpha_4=-\alpha_{46}, \; d\alpha_5=4\alpha_{56}, \; d\alpha_6 =0 $ \\
\hline
$ \mathfrak{g}_{6.56}^{\frac{4}{3}}$ & $ d\alpha_1=-\alpha_{35}-\alpha_{16}, \; d\alpha_2= -\alpha_{45}+\frac{1}{3}\alpha_{26}-\alpha_{36}, $ \\ &  $d\alpha_3=\frac{1}{3}\alpha_{36},\;  d\alpha_4=\frac{5}{3}\alpha_{46}, \; d\alpha_5=-\frac{4}{3}\alpha_{56}, \; d\alpha_6 =0 $ \\
\hline
$ \mathfrak{g}_{6.57}^{-\frac{2}{3}}$ & $ d\alpha_1=-\alpha_{35}-\alpha_{16}, \; d\alpha_2= -\alpha_{45}+\frac{4}{3}\alpha_{26}, $ \\ &  $d\alpha_3=-\frac{5}{3}\alpha_{36},\;  d\alpha_4=\frac{2}{3}\alpha_{46}-\alpha_{56}, \; d\alpha_5=\frac{2}{3}\alpha_{56}, \; d\alpha_6 =0 $ \\
\hline

\end{tabular}
}
\end{center}

\begin{center}
{\small
\begin{tabular}{|c|c|}
\hline
$ \mathfrak{g}_{6.61}^{-\frac{3}{4}}$ & $ d\alpha_1=-\alpha_{35}-2\alpha_{16}, \; d\alpha_2= -\alpha_{45}+\frac{3}{2}\alpha_{26}, $ \\ &  $d\alpha_3=-\alpha_{36}-\alpha_{56},\;  d\alpha_4=\frac{5}{2}\alpha_{46}, \; d\alpha_5=-\alpha_{56}, \; d\alpha_6 =0 $ \\
\hline
$ \mathfrak{g}_{6.63}^{-1}$ & $ d\alpha_1=-\alpha_{35}-\alpha_{16}, \; d\alpha_2= -\alpha_{45}+\alpha_{26}-\alpha_{46}, $ \\ &  $d\alpha_3=-\alpha_{36},\;  d\alpha_4=\alpha_{46}, \; d\alpha_5=0, \; d\alpha_6 =0 $ \\
\hline
$ \mathfrak{g}_{6.65}^{4l,l}$ & $ d\alpha_1=-\alpha_{35}-l\alpha_{16}, \; d\alpha_2= -\alpha_{45}-\alpha_{16}-l\alpha_{26}, $ \\ &  $d\alpha_3=3l\alpha_{36},\;  d\alpha_4=-\alpha_{36}+3l\alpha_{46}, \; d\alpha_5=-4l\alpha_{56}, \; d\alpha_6 =0 $ \\
\hline
$ \mathfrak{g}_{6.70}^{4p,p}$ & $ d\alpha_1=-\alpha_{35}-p\alpha_{16}+\alpha_{26}, \; d\alpha_2= -\alpha_{45}-\alpha_{16}-p\alpha_{26}, $ \\ &  $d\alpha_3=3p\alpha_{36}+\alpha_{46},\; d\alpha_4=-\alpha_{36}+3p\alpha_{46}, \; d\alpha_5=-4p\alpha_{56}, \; d\alpha_6 =0 $ \\
\hline
$ \mathfrak{g}_{6.71}^{-\frac{7}{4}}$ & $ d\alpha_1=-\alpha_{25}-\frac{5}{4}\alpha_{16}, \; d\alpha_2= -\alpha_{35}-\frac{1}{4}\alpha_{26}, $ \\ &  $d\alpha_3=-\alpha_{45}+\frac{3}{4}\alpha_{36},\; d\alpha_4=\frac{7}{4}\alpha_{46}, \; d\alpha_5=-\alpha_{56}, \; d\alpha_6 =0 $ \\
\hline
$ \mathfrak{g}_{6.76}^{-1}$ & $ d\alpha_1=-\alpha_{25}+\alpha_{16}, \; d\alpha_2= -\alpha_{45}, $ \\ &  $d\alpha_3=-\alpha_{24}-\alpha_{36},\; d\alpha_4=-\alpha_{46}, \; d\alpha_5=\alpha_{56}, \; d\alpha_6 =0 $ \\
\hline
$ \mathfrak{g}_{6.78}$ & $ d\alpha_1=-\alpha_{25}+\alpha_{16}, \; d\alpha_2= -\alpha_{45},\;d\alpha_3=-\alpha_{24}-\alpha_{36}-\alpha_{46}, $ \\ &  $ d\alpha_4=-\alpha_{46}, \; d\alpha_5=\alpha_{56}, \; d\alpha_6 =0 $ \\
\hline
$ \mathfrak{g}_{6.83}^{0,l}$ & $ d\alpha_1=-\alpha_{24}-\alpha_{35}, \; d\alpha_2= -l\alpha_{26}-\alpha_{36},\;d\alpha_3=-l\alpha_{36}, $ \\ &  $ d\alpha_4=l\alpha_{46}, \; d\alpha_5=\alpha_{46}+l\alpha_{56}, \; d\alpha_6 =0 $ \\
\hline
$ \mathfrak{g}_{6.84}$ & $ d\alpha_1=-\alpha_{24}-\alpha_{35}, \; d\alpha_2= -\alpha_{26},\;d\alpha_3=-\alpha_{56}, $ \\ &  $ d\alpha_4=\alpha_{46}, \; d\alpha_5=0, \; d\alpha_6 =0 $ \\
\hline
$ \mathfrak{g}_{6.88}^{0,\mu,\nu}$ & $ d\alpha_1=-\alpha_{24}-\alpha_{35}, \; d\alpha_2= -\mu\alpha_{26}+\nu\alpha_{36},\;d\alpha_3=-\nu\alpha_{26}-\mu\alpha_{36}, $ \\ &  $ d\alpha_4=\mu\alpha_{46}+\nu\alpha_{56}, \; d\alpha_5=-\nu\alpha_{46}+\mu\alpha_{56}, \; d\alpha_6 =0 $ \\
\hline
$ \mathfrak{g}_{6.89}^{0,\nu,s}$ & $ d\alpha_1=-\alpha_{24}-\alpha_{35}, \; d\alpha_2= -s\alpha_{26},\;d\alpha_3=\nu\alpha_{56}, $ \\ &  $ d\alpha_4=s\alpha_{46}, \; d\alpha_5=-\nu\alpha_{36}, \; d\alpha_6 =0 $ \\
\hline
$ \mathfrak{g}_{6.90}^{0,\nu}$ & $ d\alpha_1=-\alpha_{24}-\alpha_{35}, \; d\alpha_2= -\alpha_{46},\;d\alpha_3=\nu\alpha_{56}, $ \\{\scriptsize$ \nu\neq 1 $} &  $ d\alpha_4=-\alpha_{26}, \; d\alpha_5=-\nu\alpha_{36}, \; d\alpha_6 =0 $ \\
\hline
$ \mathfrak{g}_{6.91}$ & $ d\alpha_1=-\alpha_{24}-\alpha_{35}, \; d\alpha_2= -\alpha_{46},\;d\alpha_3=\alpha_{56}, $ \\ &  $ d\alpha_4=-\alpha_{26}, \; d\alpha_5=-\alpha_{36}, \; d\alpha_6 =0 $ \\
\hline
$ \mathfrak{g}_{6.92}^{0,\mu,\nu}$ & $ d\alpha_1=-\alpha_{24}-\alpha_{35}, \; d\alpha_2= \mu\alpha_{36},\;d\alpha_3=-\nu\alpha_{26}, $ \\ &  $ d\alpha_4=\nu\alpha_{56}, \; d\alpha_5=-\mu\alpha_{46}, \; d\alpha_6 =0 $ \\
\hline
$ \mathfrak{g}_{6.92^*}^{0}$ & $ d\alpha_1=-\alpha_{24}-\alpha_{35}, \; d\alpha_2= \alpha_{46},\;d\alpha_3=\alpha_{56}, $ \\ &  $ d\alpha_4=-\alpha_{26}, \; d\alpha_5=-\alpha_{36}, \; d\alpha_6 =0 $ \\
\hline
$ \mathfrak{g}_{6.93}^{0,\nu}$ & $ d\alpha_1=-\alpha_{24}-\alpha_{35}, \; d\alpha_2= -\alpha_{46}+\nu\alpha_{56},\;d\alpha_3=\nu\alpha_{46}, $ \\ &  $ d\alpha_4=-\alpha_{26}-\nu\alpha_{36}, \; d\alpha_5=-\nu\alpha_{26}, \; d\alpha_6 =0 $ \\
\hline
$ \mathfrak{g}_{6.94}^{-2}$ & $ d\alpha_1=-\alpha_{25}-\alpha_{34}, \; d\alpha_2=-\alpha_{35} +\alpha_{26},\;d\alpha_3=2\alpha_{36}, $ \\ &  $ d\alpha_4=-2\alpha_{46}, \; d\alpha_5=-\alpha_{56}, \; d\alpha_6 =0 $ \\
\hline
$ \mathfrak{g}_{6.101}^{a,b,c,e}$ & $ d\alpha_1=a\alpha_{15}+b\alpha_{16}, \; d\alpha_2=c\alpha_{25} +e\alpha_{26}, $ \\{\scriptsize$ a+c=-1,\; b+e=-1, $} &  $ d\alpha_3=\alpha_{36},\; d\alpha_4=\alpha_{45}, \; d\alpha_5=0, \; d\alpha_6 =0 $ \\ {\scriptsize$ ab\neq 0, \; c^2+e^2\neq 0 $} & \\
\hline
$ \mathfrak{g}_{6.102}^{-1,b,-2-b}$ & $ d\alpha_1=-\alpha_{15}+b\alpha_{16}, \; d\alpha_2=\alpha_{25} -(2+b)\alpha_{26}, $ \\ &  $ d\alpha_3=\alpha_{36},\; d\alpha_4=\alpha_{35}+\alpha_{46}, \; d\alpha_5=0, \; d\alpha_6 =0 $ \\
\hline
$ \mathfrak{g}_{6.105}^{-2,-1}$ & $ d\alpha_1=-2\alpha_{15}-\alpha_{16}, \; d\alpha_2=\alpha_{26}, \; d\alpha_3=\alpha_{35},$ \\ &  $  d\alpha_4=\alpha_{35}+\alpha_{45}, \; d\alpha_5=0, \; d\alpha_6 =0 $ \\
\hline
$ \mathfrak{g}_{6.107}^{-1,0}$ & $ d\alpha_1=-\alpha_{15}-\alpha_{26}, \; d\alpha_2=-\alpha_{25}-\alpha_{16}, \; d\alpha_3=\alpha_{35},$ \\ &  $  d\alpha_4=\alpha_{35}+\alpha_{45}, \; d\alpha_5=0, \; d\alpha_6 =0 $ \\
\hline
$ \mathfrak{g}_{6.113}^{a,b,-a,c}$ & $ d\alpha_1=a\alpha_{15}+b\alpha_{16}, \; d\alpha_2=-a\alpha_{25}+c\alpha_{26}, \; d\alpha_3=\alpha_{36},$ \\{\scriptsize$ a^2+b^2\neq 0, \;a^2+c^2\neq 0,  $} &  $  d\alpha_4=\alpha_{35}+\alpha_{46}, \; d\alpha_5=0, \; d\alpha_6 =0 $ \\ {\scriptsize$ b+c=-2 $} & \\
\hline
$ \mathfrak{g}_{6.114}^{a,-1,-\frac{a}{2}}$ & $ d\alpha_1=a\alpha_{15}-\alpha_{16}, \; d\alpha_2=\alpha_{26}, \; d\alpha_3=-\frac{a}{2}\alpha_{35}-\alpha_{45},$ \\{\scriptsize$ a\neq 0,  $} &  $  d\alpha_4=\alpha_{35}-\frac{a}{2}\alpha_{45}, \; d\alpha_5=0, \; d\alpha_6 =0 $ \\
\hline
\end{tabular}
}
\end{center}

\begin{center}
{\small
\begin{tabular}{|c|c|}
\hline
$ \mathfrak{g}_{6.115}^{-1,b,c,-c}$ & $ d\alpha_1=\alpha_{15}+c\alpha_{16}-\alpha_{26}, \; d\alpha_2=\alpha_{25}+\alpha_{16}+c\alpha_{26}, $ \\{\scriptsize$ b\neq 0,  $} &  $ d\alpha_3=-\alpha_{35}-b\alpha_{45}-c\alpha_{36},\;d\alpha_4=b\alpha_{35}-\alpha_{45}-c\alpha_{46}, \; d\alpha_5=0, \; d\alpha_6 =0 $ \\
\hline
$ \mathfrak{g}_{6.116}^{0,-1}$ & $ d\alpha_1=\alpha_{16}, \; d\alpha_2=\alpha_{15}+\alpha_{26},\;d\alpha_3=-\alpha_{45}-\alpha_{36}, $ \\ &  $ d\alpha_4=\alpha_{35}-\alpha_{46}, \; d\alpha_5=0, \; d\alpha_6 =0 $ \\
\hline
$ \mathfrak{g}_{6.118}^{0,b,-1}$ & $ d\alpha_1=-\alpha_{25}+\alpha_{16}, \; d\alpha_2=\alpha_{15}+\alpha_{26},\;d\alpha_3=-b\alpha_{45}-\alpha_{36}, $ \\ &  $ d\alpha_4=b\alpha_{35}-\alpha_{46}, \; d\alpha_5=0, \; d\alpha_6 =0 $ \\
\hline
$ \mathfrak{g}_{6.120}^{-1,-1}$ & $ d\alpha_1=-\alpha_{56}, \; d\alpha_2=-\alpha_{25}-\alpha_{26},\;d\alpha_3=\alpha_{36}, $ \\ &  $ d\alpha_4=\alpha_{45}, \; d\alpha_5=0, \; d\alpha_6 =0 $ \\
\hline
$ \mathfrak{g}_{6.125}^{0,-2}$ & $ d\alpha_1=-\alpha_{56}, \; d\alpha_2=-2\alpha_{26},\;d\alpha_3=-\alpha_{45}+\alpha_{36}, $ \\ &  $ d\alpha_4=\alpha_{35}+\alpha_{46}, \; d\alpha_5=0, \; d\alpha_6 =0 $ \\
\hline
$ \mathfrak{g}_{6.129}^{-2,-2}$ & $ d\alpha_1=-\alpha_{23}+\alpha_{15}+\alpha_{16}, \; d\alpha_2=\alpha_{25},\;d\alpha_3=\alpha_{36}, $ \\ &  $ d\alpha_4=-2\alpha_{45}-2\alpha_{46}, \; d\alpha_5=0, \; d\alpha_6 =0 $ \\
\hline
$ \mathfrak{g}_{6.135}^{0,-4}$ & $ d\alpha_1=-\alpha_{23}+2\alpha_{16}, \; d\alpha_2=\alpha_{26},\;d\alpha_3=\alpha_{25}+\alpha_{36}, $ \\ &  $ d\alpha_4=-4\alpha_{46}, \; d\alpha_5=0, \; d\alpha_6 =0 $ \\
\hline
$ \mathfrak{n}_{6.83}^{0}$ & $ d\alpha_1=-\alpha_{45}, \; d\alpha_2=-\alpha_{15}-\alpha_{36},\;d\alpha_3=-\alpha_{14}+\alpha_{26}, $ \\ &  $ d\alpha_4=\alpha_{56}, \; d\alpha_5=-\alpha_{46}, \; d\alpha_6 =0 $ \\
\hline
$ \mathfrak{n}_{6.84}^{\varepsilon}$ & $ d\alpha_1=-\alpha_{45}, \; d\alpha_2=-\alpha_{15}-\alpha_{36},\;d\alpha_3=-\alpha_{14}+\alpha_{26}-\varepsilon\alpha_{56}, $ \\{\scriptsize$ \varepsilon = \pm 1  $} &  $ d\alpha_4=\alpha_{56}, \; d\alpha_5=-\alpha_{46}, \; d\alpha_6 =0 $ \\
\hline
$ \mathfrak{n}_{6.96}^{b}$ & $ d\alpha_1=-\alpha_{24}-\alpha_{35}, \; d\alpha_2=-\alpha_{46},\;d\alpha_3=-b\alpha_{56}, $ \\ &  $ d\alpha_4=\alpha_{26}, \; d\alpha_5=-b\alpha_{36}, \; d\alpha_6 =0 $ \\
\hline
\hline
$ \mathfrak{g}_{5.7}^{p,q,r}$ & $ d\alpha_1=-\alpha_{15}, \; d\alpha_2=-p\alpha_{25},\;d\alpha_3=-q\alpha_{35}, $ \\{\scriptsize$ -1\leq r \leq q \leq p \leq 1,  $} &  $ d\alpha_4=-r\alpha_{45}, \; d\alpha_5=0 $ \\ {\scriptsize$ pqr \neq 0, \; p+q+r=-1  $} & \\
\hline
$ \mathfrak{g}_{5.8}^{-1}$ & $ d\alpha_1=-\alpha_{25}, \; d\alpha_2=0,\;d\alpha_3=-\alpha_{35}, \; d\alpha_4=\alpha_{45}, \; d\alpha_5=0$ \\
\hline
$ \mathfrak{g}_{5.9}^{p,-2-p}$ & $ d\alpha_1=-\alpha_{15}-\alpha_{25}, \; d\alpha_2=-\alpha_{25},\;d\alpha_3=-p\alpha_{35}, $ \\{\scriptsize$ p\geq -1,  $} &  $ d\alpha_4=-(p+2)\alpha_{45}, \; d\alpha_5=0 $ \\
\hline
$ \mathfrak{g}_{5.11}^{-3}$ & $ d\alpha_1=-\alpha_{15}-\alpha_{25}, \; d\alpha_2=-\alpha_{25}-\alpha_{35},\;d\alpha_3=-\alpha_{35}, $ \\ &  $ d\alpha_4=3\alpha_{45}, \; d\alpha_5=0 $ \\
\hline
$ \mathfrak{g}_{5.13}^{-1-2q,q,r}$ & $ d\alpha_1=-\alpha_{15}, \; d\alpha_2=(1+2q)\alpha_{25},\;d\alpha_3=-q\alpha_{35}-r\alpha_{45}, $ \\{\scriptsize$ -1\leq q \leq 0,  $} &  $ d\alpha_4=r\alpha_{35}-q\alpha_{45}, \; d\alpha_5=0 $ \\ {\scriptsize$ q \neq -\frac{1}{2}, \; r\neq 0  $} & \\
\hline
$ \mathfrak{g}_{5.14}^{0}$ & $ d\alpha_1=-\alpha_{25}, \; d\alpha_2=0,\;d\alpha_3=-\alpha_{45},\; d\alpha_4=\alpha_{35}, \; d\alpha_5=0 $ \\
\hline
$ \mathfrak{g}_{5.15}^{-1}$ & $ d\alpha_1=-\alpha_{15}-\alpha_{25}, \; d\alpha_2=-\alpha_{25},\;d\alpha_3=\alpha_{35}-\alpha_{45}, $ \\ &  $ d\alpha_4=\alpha_{45}, \; d\alpha_5=0 $ \\
\hline
$ \mathfrak{g}_{5.16}^{-1,q}$ & $ d\alpha_1=-\alpha_{15}-\alpha_{25}, \; d\alpha_2=-\alpha_{25},\;d\alpha_3=\alpha_{35}-q\alpha_{45}, $ \\{\scriptsize$ q \neq 0  $} &  $ d\alpha_4=q\alpha_{35}+\alpha_{45}, \; d\alpha_5=0 $ \\
\hline
$ \mathfrak{g}_{5.17}^{p,-p,r}$ & $ d\alpha_1=-p\alpha_{15}-\alpha_{25}, \; d\alpha_2=\alpha_{15}-p\alpha_{25},\;d\alpha_3=p\alpha_{35}-r\alpha_{45}, $ \\{\scriptsize$ r \neq 0  $} &  $ d\alpha_4=r\alpha_{35}+p\alpha_{45}, \; d\alpha_5=0 $ \\
\hline
$ \mathfrak{g}_{5.18}^{0}$ & $ d\alpha_1=-\alpha_{25}-\alpha_{35}, \; d\alpha_2=\alpha_{15}-\alpha_{45},\;d\alpha_3=-\alpha_{45}, $ \\ &  $ d\alpha_4=\alpha_{35}, \; d\alpha_5=0$ \\
\hline
$ \mathfrak{g}_{5.19}^{p,-2p-2}$ & $ d\alpha_1=-\alpha_{23}-(p+1)\alpha_{15}, \; d\alpha_2=-\alpha_{25},\;d\alpha_3=-p\alpha_{35}, $ \\{\scriptsize$ p \neq -1  $} &  $ d\alpha_4=(2p+2)\alpha_{45}, \; d\alpha_5=0$ \\
\hline
$ \mathfrak{g}_{5.20}^{-1}$ & $ d\alpha_1=-\alpha_{23}-\alpha_{45}, \; d\alpha_2=-\alpha_{25},\;d\alpha_3=\alpha_{35}, \; d\alpha_4=0, \; d\alpha_5=0 $ \\
\hline
$ \mathfrak{g}_{5.23}^{-4}$ & $ d\alpha_1=-\alpha_{23}-2\alpha_{15}, \; d\alpha_2=-\alpha_{25},\;d\alpha_3=-\alpha_{25}-\alpha_{35}, $ \\ &  $ d\alpha_4=4\alpha_{45}, \; d\alpha_5=0$ \\
\hline
$ \mathfrak{g}_{5.25}^{p,4p}$ & $ d\alpha_1=-\alpha_{23}-2p\alpha_{15}, \; d\alpha_2=-p\alpha_{25}+\alpha_{35},\;d\alpha_3=-\alpha_{25}-p\alpha_{35}, $ \\ {\scriptsize$ p \neq 0  $}  &  $ d\alpha_4=4p\alpha_{45}, \; d\alpha_5=0 $ \\
\hline
$ \mathfrak{g}_{5.26}^{0,\varepsilon}$ & $ d\alpha_1=-\alpha_{23}-\varepsilon\alpha_{45}, \; d\alpha_2=\alpha_{35},\;d\alpha_3=-\alpha_{25}, $ \\ {\scriptsize$ \varepsilon = \pm 1  $}  &  $ d\alpha_4=0, \; d\alpha_5=0$ \\
\hline
$ \mathfrak{g}_{5.28}^{-\frac{3}{2}}$ & $ d\alpha_1=-\alpha_{23}+\frac{1}{2}\alpha_{15}, \; d\alpha_2=\frac{3}{2}\alpha_{25},\;d\alpha_3=-\alpha_{35}, $ \\   &  $ d\alpha_4=-\alpha_{35}-\alpha_{45}, \; d\alpha_5=0$ \\
\hline
\end{tabular}
}
\end{center}

\begin{center}
{\small
\begin{tabular}{|c|c|}
\hline
$ \mathfrak{g}_{5.30}^{-\frac{4}{3}}$ & $ d\alpha_1=-\alpha_{24}-\frac{2}{3}\alpha_{15}, \; d\alpha_2=-\alpha_{34}+\frac{1}{3}\alpha_{25},\;d\alpha_3=\frac{4}{3}\alpha_{35}, $ \\   &  $ d\alpha_4=-\alpha_{45}, \; d\alpha_5=0$ \\
\hline
$ \mathfrak{g}_{5.33}^{-1,-1}$ & $ d\alpha_1=-\alpha_{14}, \; d\alpha_2=-\alpha_{25},\;d\alpha_3=\alpha_{34}+\alpha_{35}, $ \\   &  $ d\alpha_4=0, \; d\alpha_5=0$ \\
\hline
$ \mathfrak{g}_{5.35}^{-2,0}$ & $ d\alpha_1=2\alpha_{14}, \; d\alpha_2=-\alpha_{24}-\alpha_{35},\;d\alpha_3=-\alpha_{34}+\alpha_{25}, $ \\   &  $ d\alpha_4=0, \; d\alpha_5=0 $ \\
\hline
\hline
$ \mathfrak{g}_{4.2}^{-2}$ & $ d\alpha_1=2\alpha_{14}, \; d\alpha_2=-\alpha_{24}-\alpha_{34}, $ \\ &  $ d\alpha_3=-\alpha_{34},\;d\alpha_4=0 $ \\ 
\hline
$ \mathfrak{g}_{4.5}^{p,-p-1}$ & $ d\alpha_1=-\alpha_{14}, \; d\alpha_2=-p\alpha_{24}, $ \\{\scriptsize$ -1 \leq p \leq -\frac{1}{2},  $} &  $d\alpha_3=(p+1)\alpha_{34},\; d\alpha_4=0 $ \\ 
\hline
$ \mathfrak{g}_{4.6}^{-2p,p}$ & $ d\alpha_1=2p\alpha_{14}, \; d\alpha_2=-p\alpha_{24}-\alpha_{34}, $ \\{\scriptsize$  p> 0  $} &  $d\alpha_3=\alpha_{24}-p\alpha_{34},\; d\alpha_4=0 $ \\ 
\hline
$ \mathfrak{g}_{4.8}^{-1}$ & $ d\alpha_1=-\alpha_{23}, \; d\alpha_2=-\alpha_{24},\;d\alpha_3=\alpha_{34}, \; d\alpha_4=0 $ \\
\hline
$ \mathfrak{g}_{4.9}^{0}$ & $ d\alpha_1=-\alpha_{23}, \; d\alpha_2=-\alpha_{34},\;d\alpha_3=\alpha_{24},\; d\alpha_4=0 $ \\
\hline
\hline
$ \mathfrak{g}_{3.1}$ & $ d\alpha_1=-\alpha_{23}, \; d\alpha_2=0,\;d\alpha_3=0, $ \\ {\scriptsize nilpotent}  &   \\
\hline
$ \mathfrak{g}_{3.4}^{-1}$ & $ d\alpha_1=-\alpha_{13}, \; d\alpha_2=\alpha_{23},\;d\alpha_3=0, $ \\ 
\hline
$ \mathfrak{g}_{3.5}^{0}$ & $ d\alpha_1=-\alpha_{23}, \; d\alpha_2=\alpha_{13},\;d\alpha_3=0, $ \\ 
\hline
\end{tabular}
}
\end{center}

\newpage

\section*{Table 3: { \rm Symplectic structures on 6-dimensional solvable unimodular Lie algebras}} 

\begin{center}
{\small
\begin{tabular}{|c|c|c|}
\hline
Lie algebra & Symplectic form & Conditions on $\omega_{i,j}$ \\
\hline 
$\mathfrak{g}_{6.3}^{0,-1}$ & $ \omega = \omega_{1,6}\alpha_{16}+\omega_{2,3}\alpha_{23}+\omega_{2,6}\alpha_{26}+\omega_{3,6}\alpha_{36} $ & $ \omega_{1,6}\omega_{2,3}\omega_{4,5} \neq 0$ \\
& $+\omega_{4,5}\alpha_{45}+\omega_{4,6}\alpha_{46}+\omega_{5,6}\alpha_{56}$ & \\
\hline 
$\mathfrak{g}_{6.10}^{0,0}$ & $\omega = \omega_{1,6}\alpha_{16}+\omega_{2,3}\alpha_{23}+\omega_{2,6}\alpha_{26}+\omega_{3,6}\alpha_{36}$ & $ \omega_{1,6}\omega_{2,3}\omega_{4,5} \neq 0$ \\
& $+\omega_{4,5}\alpha_{45}+\omega_{4,6}\alpha_{46}+\omega_{5,6}\alpha_{56} $ & \\
\hline 
$\mathfrak{g}_{6.13}^{\frac{1}{2},-1,0}$ & $ \omega = \omega_{1,2}\alpha_{12}+\omega_{2,3}(-\frac{1}{2}\alpha_{16}+\alpha_{23})+\omega_{2,6}\alpha_{26}$ & $ \omega_{1,2}\omega_{3,4}\omega_{5,6} \neq 0$ \\
& $+\omega_{3,4}\alpha_{34}+\omega_{3,6}\alpha_{36}+\omega_{4,6}\alpha_{46}+\omega_{5,6}\alpha_{56} $ & \\
\hline 
$\mathfrak{g}_{6.13}^{-1,\frac{1}{2},0}$ & $ \omega = \omega_{1,3}\alpha_{13}+\omega_{2,3}(-\frac{1}{2}\alpha_{16}+\alpha_{23})+\omega_{2,4}\alpha_{24}+\omega_{2,6}\alpha_{26}$ & $\omega_{1,3}\omega_{2,4}\omega_{5,6} \neq 0$\\
& $+\omega_{3,6}\alpha_{36}+\omega_{4,6}\alpha_{46}+\omega_{5,6}\alpha_{56} $ & \\
\hline 
$\mathfrak{g}_{6.15}^{-1}$ & $ \omega =\omega_{1,6}\alpha_{16}+(\frac{\omega_{1,6}+\omega_{3,4}}{2})\alpha_{25}+\omega_{3,4}\alpha_{34}$ & $ \omega_{1,6}\omega_{3,4} \neq 0, $\\
& $+\omega_{3,6}\alpha_{36}+\omega_{4,6}\alpha_{46}+\omega_{5,6}\alpha_{56} $ & $\omega_{1,6}\neq -\omega_{3,4} $ \\
\hline 
$\mathfrak{g}_{6.18}^{-1,-1}$ & $ \omega = \omega_{1,6}(\alpha_{16}+\alpha_{24})+\omega_{2,6}\alpha_{26}+\omega_{3,5}\alpha_{35}$ & $ \omega_{1,6}\omega_{3,5} \neq 0$\\
& $+\omega_{3,6}\alpha_{36}+\omega_{4,6}\alpha_{46}+\omega_{5,6}\alpha_{56} $ & \\
\hline
$\mathfrak{g}_{6.21}^{0}$ & $ \omega = \omega_{1,2}\alpha_{12}+\omega_{2,3}\alpha_{23}+\omega_{2,6}\alpha_{26}+\omega_{3,6}\alpha_{36}$ & $ \omega_{1,2}\omega_{3,6}\omega_{4,5} \neq 0$\\
& $+\omega_{4,5}\alpha_{45}+\omega_{4,6}\alpha_{46}+\omega_{5,6}\alpha_{56} $ & \\
\hline 
$\mathfrak{g}_{6.23}^{0,0,\varepsilon}$ & $ \omega =\omega_{1,2}(\alpha_{12}+\varepsilon \alpha_{35})+\omega_{1,6}(\alpha_{16}+ \alpha_{24})+\omega_{2,3}\alpha_{23}$ & $ \omega_{1,2}\neq 0, $\\
& $+\omega_{2,5}\alpha_{25}+\omega_{2,6}\alpha_{26}+\omega_{3,6}\alpha_{36}+\omega_{4,6}\alpha_{46}+\omega_{5,6}\alpha_{56} $ & $\omega_{1,6}^2 + \omega_{1,2} \omega_{4,6}\neq 0$ \\
\hline 
$\mathfrak{g}_{6.29}^{0,0,\varepsilon \neq 0}$ & $ \omega =\omega_{1,3}(\alpha_{13}+\varepsilon \alpha_{45})+\omega_{1,6}(\alpha_{16}+ \alpha_{24})+\omega_{2,3}\alpha_{23}+\omega_{2,6}\alpha_{26}$ & $ \omega_{1,3}\neq 0, $\\
& $+\omega_{3,4}\alpha_{34}+\omega_{3,6}\alpha_{36}+\omega_{4,6}\alpha_{46}+\omega_{5,6}\alpha_{56} $ & $ (\omega_{1,6},\omega_{2,6}) \neq (0,0), $ \\
&& $\omega_{5,6}\omega_{1,6} \neq \varepsilon(\omega_{2,3}\omega_{1,6}-$ \\
&& $\omega_{2,6}\omega_{1,3})$ \\
\hline
$\mathfrak{g}_{6.29}^{0,0,0}$ & $ \omega =\omega_{1,2}\alpha_{12}+\omega_{1,3}\alpha_{13}+\omega_{1,6}(\alpha_{16}+ \alpha_{24})+\omega_{2,3}\alpha_{23}+\omega_{2,6}\alpha_{26}$ & $ \omega_{5,6}\neq 0, $\\
& $+\omega_{3,4}\alpha_{34}+\omega_{3,6}\alpha_{36}+\omega_{4,6}\alpha_{46}+\omega_{5,6}\alpha_{56} $ & $\omega_{1,3}\omega_{1,6}-\omega_{1,2}\omega_{3,4}\neq 0 $\\
\hline 
$\mathfrak{g}_{6.36}^{0,0}$ & $ \omega =\omega_{1,2}\alpha_{12}+\omega_{2,3}\alpha_{23}+\omega_{2,6}\alpha_{26}+\omega_{3,6}\alpha_{36} $ & $  \omega_{1,2}\omega_{3,6}\omega_{4,5}\neq 0 $\\
& $+\omega_{4,5}\alpha_{45}+\omega_{4,6}\alpha_{46}+\omega_{5,6}\alpha_{56}$ & \\
\hline 
$\mathfrak{g}_{6.38}^{0}$ & $ \omega=\omega_{1,6}(2\alpha_{16}+\alpha_{25}-\alpha_{34})+  $ & $ \omega_{1,6}\neq 0$\\
& $\omega_{2,3}\alpha_{23}+\omega_{2,4}(\alpha_{24}+\alpha_{35})+\omega_{2,6}\alpha_{26}$ &\\
& $+\omega_{3,6}\alpha_{36}+\omega_{4,6}\alpha_{46}+\omega_{5,6}\alpha_{56}$ &\\
\hline 
$\mathfrak{g}_{6.54}^{0,-1}$ & $ \omega=\omega_{1,4}(\alpha_{14}+\alpha_{23})+\omega_{1,6}(\alpha_{16}+\alpha_{35})+\omega_{2,6}(\alpha_{26}- $ & $\omega_{1,4}\omega_{5,6}\neq 0$\\
& $\alpha_{45})+\omega_{3,4}\alpha_{34}+\omega_{3,6}\alpha_{36}+\omega_{4,6}\alpha_{46}+\omega_{5,6}\alpha_{56}$ & \\
\hline
$\mathfrak{g}_{6.70}^{0,0}$ & $ \omega=\omega_{1,3}(\alpha_{13}+\alpha_{24})+\omega_{1,6}(\alpha_{16}+\alpha_{45})+\omega_{3,4}\alpha_{34}+ $ & $\omega_{1,3}\omega_{1,6}\omega_{3,5}\omega_{5,6}\neq 0, $\\
& $\omega_{3,5}\alpha_{35}+\omega_{3,6}\alpha_{36}+\omega_{4,6}\alpha_{46}+\omega_{5,6}\alpha_{56}$ & $\omega_{1,3}\omega_{5,6} + \omega_{1,6}\omega_{3,5} \neq 0 $ \\
\hline 
$\mathfrak{g}_{6.78}$ & $\omega=\omega_{1,4}(\alpha_{14}+\alpha_{35}-\alpha_{26})+\omega_{1,6}(\alpha_{16}-\alpha_{25}) $ & $ \omega_{1,4} \neq 0$\\
& $+\omega_{2,4}(\alpha_{24}+\alpha_{36})+\omega_{4,5}\alpha_{45}+\omega_{4,6}\alpha_{46}+\omega_{5,6}\alpha_{56}$ & \\
\hline 
$\mathfrak{g}_{6.118}^{0,\pm 1,-1}$ & $  \omega=\omega_{1,3}(\alpha_{13}\pm\alpha_{24})+\omega_{1,4}(\alpha_{14}\mp\alpha_{23})+\omega_{1,5}(\alpha_{15}+$ & $ \omega_{1,4}\omega_{1,6}\omega_{5,6}\neq 0, $\\
 & $\alpha_{26})+\omega_{1,6}(\alpha_{16}+\alpha_{25})+\omega_{3,5}(\pm\alpha_{35}-\alpha_{46})+ $ & $2\omega_{1,4}\omega_{1,6}\omega_{3,6} + \omega_{1,3}\omega_{3,5}\pm $ \\
 & $\omega_{3,6}(\alpha_{36}\pm \alpha_{45})+\omega_{5,6}\alpha_{56}$ & $\omega_{1,3}^2\omega_{5,6}\pm \omega_{1,4}^2\omega_{5,6}\neq 0$ \\
\hline 
$\mathfrak{n}_{6.84}$ & $ \omega=\omega_{1,4}(\alpha_{14}-\alpha_{26})+\omega_{1,5}(\alpha_{15}+\alpha_{36})+\omega_{1,6}(-\varepsilon\alpha_{16}+\alpha_{25}+$ & $\omega_{1,6}\neq 0$\\
& $\alpha_{34})+\omega_{4,5}\alpha_{45}+\omega_{4,6}\alpha_{46}+\omega_{5,6}\alpha_{56} $ & \\
\hline 
$\mathfrak{g}_{5.7}^{p,-p,-1}\oplus \mathbb{R}$ & $ \omega =\omega_{1,4}\alpha_{14}+\omega_{1,5}\alpha_{15}+\omega_{2,3}\alpha_{23}+\omega_{2,5}\alpha_{25}+$ & $ \omega_{1,4}\omega_{2,3}\omega_{5,6}\neq 0$\\
& $\omega_{3,5}\alpha_{35}+\omega_{4,5}\alpha_{45}+\omega_{5,6}\alpha_{56} $ & \\
\hline 
$\mathfrak{g}_{5.7}^{1,-1,-1}\oplus \mathbb{R}$ & $ \omega = \omega_{1,3}\alpha_{13}+\omega_{1,4}\alpha_{14}+\omega_{1,5}\alpha_{15}+\omega_{2,3}\alpha_{23}+ $ & $ \omega_{1,4}\omega_{2,3}-\omega_{1,3}\omega_{2,4}\neq 0,$\\
& $\omega_{2,4}\alpha_{24}+\omega_{2,5}\alpha_{25}+\omega_{3,5}\alpha_{35}+\omega_{4,5}\alpha_{45}+\omega_{5,6}\alpha_{56}$ & $ \omega_{5,6} \neq 0$ \\
\hline 
$\mathfrak{g}_{5.8}^{-1}\oplus \mathbb{R}$ & $ \omega =\omega_{1,2}\alpha_{12}+\omega_{1,5}\alpha_{15}+\omega_{2,5}\alpha_{25}+\omega_{2,6}\alpha_{26} $ & $ \omega_{3,4}\neq 0, $\\
& $+\omega_{3,4}\alpha_{34}+\omega_{3,5}\alpha_{35}+\omega_{4,5}\alpha_{45}+\omega_{5,6}\alpha_{56}$ & $\omega_{1,2}\omega_{5,6}-\omega_{1,5}\omega_{2,6}\neq 0 $ \\
\hline

\end{tabular} 
}
\end{center}

\begin{center}
{\small
\begin{tabular}{|c|c|c|}
\hline
Lie algebra & Symplectic form & Conditions on $\omega_{i,j}$ \\
\hline
$\mathfrak{g}_{5.13}^{-1,0,r}\oplus \mathbb{R}$ & $ \omega =\omega_{1,2}\alpha_{12}+\omega_{1,5}\alpha_{15}+\omega_{2,5}\alpha_{25}+\omega_{3,4}\alpha_{34}$ & $ \omega_{1,2}\omega_{3,4}\omega_{5,6}\neq 0$\\
& $+\omega_{3,5}\alpha_{35}+\omega_{4,5}\alpha_{45}+\omega_{5,6}\alpha_{56} $ & \\
\hline 
$\mathfrak{g}_{5.14}^{0}\oplus \mathbb{R}$ & $\omega =\omega_{1,2}\alpha_{12}+\omega_{1,5}\alpha_{15}+\omega_{2,5}\alpha_{25}+\omega_{2,6}\alpha_{26} $ & $\omega_{1,5}\omega_{2,6}\omega_{3,4}\omega_{5,6}\neq 0, $\\
& $+\omega_{3,4}\alpha_{34}+\omega_{3,5}\alpha_{35}+\omega_{4,5}\alpha_{45}+\omega_{5,6}\alpha_{56}$ & $\omega_{1,2}\omega_{5,6}-\omega_{1,5}\omega_{2,6}\neq 0 $ \\
\hline 
$\mathfrak{g}_{5.15}^{-1}\oplus \mathbb{R}$ & $ \omega =\omega_{1,4}(\alpha_{14}-\alpha_{23})+\omega_{1,5}\alpha_{15}+\omega_{2,4}\alpha_{24} $ & $ \omega_{1,4}\omega_{5,6}\neq 0$\\
& $+\omega_{2,5}\alpha_{25}+\omega_{3,5}\alpha_{35}+\omega_{4,5}\alpha_{45}+\omega_{5,6}\alpha_{56}$ & \\
\hline 
$\mathfrak{g}_{5.17}^{0,0,r}\oplus \mathbb{R}$ & $ \omega =\omega_{1,2}\alpha_{12}+\omega_{1,5}\alpha_{15}+\omega_{2,5}\alpha_{25}+\omega_{3,4}\alpha_{34}+$ & $ \omega_{1,2}\omega_{3,4}\omega_{5,6}\neq 0$\\
& $\omega_{3,5}\alpha_{35}+\omega_{4,5}\alpha_{45}+\omega_{5,6}\alpha_{56} $ & \\
\hline 
$\mathfrak{g}_{5.17}^{p,-p,\pm 1}\oplus \mathbb{R}$ & $\omega =\omega_{1,3}(\pm\alpha_{13}+\alpha_{24})+\omega_{1,4}(\mp\alpha_{14}+\alpha_{23})+ $ & $\omega_{1,3}\mp\omega_{1,4}\neq 0, \; \omega_{5,6} \neq 0$\\
& $\omega_{1,5}\alpha_{15}+\omega_{2,5}\alpha_{25}+\omega_{3,5}\alpha_{35}+\omega_{4,5}\alpha_{45}+\omega_{5,6}\alpha_{56}$ & \\
\hline
$\mathfrak{g}_{5.17}^{0,0,\pm 1}\oplus \mathbb{R}$ & $\omega =\omega_{1,2}\alpha_{12}+\omega_{1,3}(\pm\alpha_{13}+\alpha_{24})+$ & $ \omega_{1,2}\omega_{3,4}-\omega_{1,3}^2-\omega_{1,4}^2\neq 0, $\\
& $\omega_{1,4}(\mp\alpha_{14}+\alpha_{23})+\omega_{1,5}\alpha_{15}+\omega_{2,5}\alpha_{25}+ $ & $\omega_{3,4} \omega_{5,6} \neq 0$ \\
& $\omega_{3,4}\alpha_{34}+\omega_{3,5}\alpha_{35}+\omega_{4,5}\alpha_{45}+\omega_{5,6}\alpha_{56}$ & \\
\hline 
$\mathfrak{g}_{5.18}^{0}\oplus \mathbb{R}$
& $\omega =\omega_{1,3}(\alpha_{13}+\alpha_{24})+\omega_{1,5}\alpha_{15}+\omega_{2,5}\alpha_{25}+$ & $\omega_{1,3} \omega_{5,6} \neq 0$\\
& $\omega_{3,4}\alpha_{34}+\omega_{3,5}\alpha_{35}+\omega_{4,5}\alpha_{45}+\omega_{5,6}\alpha_{56} $ & \\
\hline 
$\mathfrak{g}_{5.19}^{-2,2}\oplus \mathbb{R}$ & $\omega =\omega_{1,2}\alpha_{12}+\omega_{1,5}(\alpha_{15}-\alpha_{23})+\omega_{2,5}\alpha_{25}+$ & $ \omega_{1,2}\omega_{3,4} \omega_{5,6} \neq 0$\\
& $\omega_{3,4}\alpha_{34}+\omega_{3,5}\alpha_{35}+\omega_{4,5}\alpha_{45}+\omega_{5,6}\alpha_{56} $ & \\
\hline 
$\mathfrak{g}_{5.19}^{-\frac{1}{2},-1}\oplus \mathbb{R}$ & $\omega =\omega_{1,3}\alpha_{13}+\omega_{1,5}(\alpha_{15}-\alpha_{23})+\omega_{2,4}\alpha_{24}+ $ & $ \omega_{1,3}\omega_{2,4} \omega_{5,6} \neq 0$\\
& $\omega_{2,5}\alpha_{25}+\omega_{3,5}\alpha_{35}+\omega_{4,5}\alpha_{45}+\omega_{5,6}\alpha_{56}$ & \\
\hline 
$\mathfrak{g}_{3.4}^{-1}\oplus 3\mathbb{R} $ & $\omega =\omega_{1,2}\alpha_{12}+\omega_{1,3}\alpha_{13}+\omega_{2,3}\alpha_{23}+\omega_{3,4}\alpha_{34}+$ & $\omega_{1,2}\omega_{5,6}\omega_{3,4}-\omega_{1,2}\omega_{4,6}\omega_{3,5}+$\\
$\mathfrak{g}_{3.5}^{0}\oplus 3\mathbb{R}$ & $\omega_{3,5}\alpha_{35}+\omega_{3,6}\alpha_{36}+\omega_{4,5}\alpha_{45}+\omega_{4,6}\alpha_{46}+\omega_{5,6}\alpha_{56} $ & $\omega_{1,2}\omega_{4, 5}\omega_{3,6}  \neq 0$ \\
\hline 
$\mathfrak{g}_{3.1}\oplus \mathfrak{g}_{3.4}^{-1}$ & $\omega =\omega_{1,2}\alpha_{12}+\omega_{1,3}\alpha_{13}+\omega_{2,3}\alpha_{23}+\omega_{2,6}\alpha_{26}+$ & $\omega_{4,5}\neq 0,  $\\
$ \mathfrak{g}_{3.1}\oplus \mathfrak{g}_{3.5}^{0}$ & $\omega_{3,6}\alpha_{36}+\omega_{4,5}\alpha_{45}+\omega_{4,6}\alpha_{46}+\omega_{5,6}\alpha_{56} $ & $\omega_{3, 6}\omega_{1, 2}-\omega_{2, 6}\omega_{1, 3} \neq 0$ \\
\hline 
$\mathfrak{g}_{3.4}^{-1}\oplus \mathfrak{g}_{3.4}^{-1}$ & $\omega =\omega_{1,2}\alpha_{12}+\omega_{1,3}\alpha_{13}+\omega_{2,3}\alpha_{23}+$ & $ \omega_{1,2}\omega_{3,6}\omega_{4,5} \neq 0 $\\
$ \mathfrak{g}_{3.4}^{-1}\oplus \mathfrak{g}_{3.5}^{0} $ & $\omega_{3,6}\alpha_{36}+\omega_{4,5}\alpha_{45}+\omega_{4,6}\alpha_{46}+\omega_{5,6}\alpha_{56} $ & \\
$ \mathfrak{g}_{3.5}^{0}\oplus \mathfrak{g}_{3.5}^{0} $ & & \\
\hline 
\end{tabular} 
}
\end{center}


\bigskip 


\begin{thebibliography}{biblio}
\bibitem{bock} C. Bock, \textquotedblleft\textit{On Low-Dimensional Solvmanifolds}\textquotedblright, Ph.D Thesis, Erlanghen University (2010).
\bibitem{cf} S. Console, A. Fino, \textquotedblleft\textit{On the de Rham cohomology of solvmanifolds}\textquotedblright, to appear in Annali della Scuola Normale Superiore di Pisa (2011).
\bibitem{FOT} Y. Felix, J. Oprea, D. Tanr\'e, Algebraic models in geometry, Oxford Graduate Texts in Mathematics, 17. Oxford University Press, Oxford, 2008.
\bibitem{FS} M. Freibert, F. Schulte-Hengesbach, \textquotedblleft\textit{Half-flat structures on indecomposable Lie groups}\textquotedblright, arXiv:1110.1512v1 [math.DG] (2011).
\bibitem{guan2} Z. D. Guan, \textquotedblleft\textit{Modification and the cohomology groups of compact solvmanifolds}\textquotedblright, Electron. Res. Announc. Amer. Math. Soc. 13 (2007), 74--81. 
\bibitem{guan} Z.-D. Guan,\textquotedblleft\textit{Toward a Classification of Compact Nilmanifolds with Symplectic Structures}\textquotedblright, International Mathematics Research Notices, Vol. 2010, No. 22, pp. 4377--4384, 2010.
\bibitem{hattori} A. Hattori, \textquotedblleft\textit{Spectral sequence in the de Rham cohomology of fibre bundles}\textquotedblright, J. Fac. Sci. Univ. Tokio \textbf{8 (Sect. 1)}, (1960), 289-331.
\bibitem{mad-sw1} T. B. Madsen and A. F. Swann, \textquotedblleft\textit{Homogeneous spaces, multi-moment maps and (2,3)-trivial algebras}\textquotedblright, 2010, IMADA preprint, CP3-ORIGINS: 2010-52, eprint arXiv:1012.0402 [math.DG]. Proceedings of the XIXth International Fall Workshop on Geometry and Physics, Porto, September 6-9, 2010, AIP Conference Proceedings, to appear.
\bibitem{mad-sw2} T. B. Madsen and A. F. Swann, \textquotedblleft\textit{Multi-moment maps}\textquotedblright, 2010, IMADA preprint, CP3-ORIGINS: 2010-53.
\bibitem{Ma} A. Mal'\v{c}ev,  \textquotedblleft\textit{On a class of homogeneous spaces}\textquotedblright, 
Amer. Math. Soc. Trasl.,  \textbf{39} (1951).
\bibitem{mostow2}  G.  Mostow, \textquotedblleft\textit{Cohomology of topological groups and solvmanifolds}\textquotedblright, Ann. of Math. (2) {\textbf 73} (1961), 20--48.
\bibitem{muba} G. M. Mubarakzyanov, \textquotedblleft\textit{On solvable Lie algebras}\textquotedblright, Izv. Vys˘s. U˘cehn. Zaved. Matematika. 32, (1963), 104-116 (Russian).
\bibitem{nomizu}  K. Nomizu, \textquotedblleft\textit{On the cohomology of compact homogeneous space of nilpotent Lie group}\textquotedblright, Ann. of Math. (2) \textbf{59} (1954), 531-538.
\bibitem{rag} M.S. Raghunathan,  \textquotedblleft\textit{Discrete Subgroups of Lie Groups}\textquotedblright, Springer-Verlag, Berlin, Heidelberg, New York, (1972).
\bibitem{salamon} 
 S. Salamon,  \textquotedblleft\textit{Complex structures on nilpotent Lie algebras}\textquotedblright, J. Pure Appl. Algebra {\bf 157} (2001), no. 2-3, 311--333. 
\bibitem{shaba} A. Shabanskaya, \textquotedblleft\textit{Classification of Six Dimensional Solvable Indecomposable Lie Algebras with a codimension one nilradical over $\mathbb{R}$}\textquotedblright, Ph.D Thesis, University of Toledo (2011).
\bibitem{thurston}  W. Thurston, \textquotedblleft\textit{Some simple examples of symplectic manifolds.}\textquotedblright, Proceedings of the American Mathematical Society 55 (1976): 467–468.
\bibitem{tsengyau} L.S. Tseng and S.T. Yau,\textquotedblleft\textit{Cohomology and Hodge theory on symplectic manifolds: I}\textquotedblright, arXiv:0909.5418v1[math.SG] (2009).
\end{thebibliography}
\end{document}